\documentclass[11pt]{amsart}
\usepackage{amsmath, amsthm, amscd, amsfonts, amssymb, graphicx, color}
\usepackage[bookmarksnumbered, colorlinks, plainpages]{hyperref}
\usepackage{plain}
\usepackage{graphicx} 
\usepackage{xcolor}

\usepackage{tikz}
\usetikzlibrary{calc}

\usepackage{float}

\usepackage{amsmath,amssymb,amsthm}

\usepackage{thmtools}
\usepackage{cleveref}

\newtheorem{theorem}{Theorem}[section]
\newtheorem{definition}[theorem]{Definition}

\newtheorem{lemma}[theorem]{Lemma}

\newtheorem{proposition}[theorem]{Proposition}
\newtheorem*{claim}{Claim}

\newtheorem{remark}[theorem]{Remark}

\newcommand{\cH}{\mathcal{H}}
\renewcommand{\L}{\mathbb{L}}

\def\N{\mathbb{N}}
\def\Z{\mathbb{Z}}
\def\Q{\mathbb{Q}}
\def\R{\mathbb{R}}

\DeclareMathOperator{\diam}{diam}

\title[Dimension of irrationals with Lagrange value exactly 3]{Generalized Hausdorff measure of irrationals with Lagrange value exactly 3}

\author[C. G. Moreira]{Carlos Gustavo Moreira}
\address{Carlos Gustavo Moreira: SUSTech International Center for Mathematics, Shenzhen, Guangdong, P. R. China and Instituto Nacional de Matemática Pura e Aplicada, Rio de Janeiro, Brasil.}
\email{gugu@impa.br}

\author[H. Erazo]{Harold Erazo}
\address{Harold Erazo: Instituto Nacional de Matemática Pura e Aplicada, Rio de Janeiro, Brasil.}
\email{harold.eraz@gmail.com}

\author{Nicolas Angelini}
\address{Instituto de Matemática Aplicada San Luis (IMASL, CONICET), San Luis, Argentina\\
and Departamento de Matemática, Facultad de Ciencias Físico Matemáticas y Naturales, Universidad Nacional de San Luis, Argentina.}
\email{nicolas.angelini.2015@gmail.com}


\begin{document}

\begin{abstract}
We study the generalized Hausdorff measure of some natural subsets of $k^{-1}(3)$, where $k^{-1}(3)$ consists of the real numbers $x$ for which  $\left| x-\frac{p}{q} \right|<\frac{1}{(3+\varepsilon)q^2}$ has infinitely many rational solutions $\frac{p}{q}$ for any $\varepsilon<0$ but only finitely many for any $\varepsilon>0$. It is well known that $k^{-1}(3)$ is an uncountable set with Hausdorff dimension zero. Given any dimension function $h$, we determine the exact ``cut point" at which the generalized Hausdorff measure $\mathcal{H}^h(k^{-1}(3))$ drops from infinity to zero. In particular we show that such a measure is always zero or not $\sigma$--finite, and, as an application, we can classify topologically $k^{-1}(3)$. Moreover, we show that the subset of attainable elements of $k^{-1}(3)$ has the same generalized Hausdorff measure as $k^{-1}(3)$, but the subset of non--attainable elements of $k^{-1}(3)$ has a ``strictly smaller" generalized Hausdorff measure.
\end{abstract}

\maketitle
\section{Introduction}

The Hausdorff dimension provides a powerful tool for quantifying the fractal structure of sets in metric spaces, but it is often too coarse to distinguish between certain extremely thin sets, all of which may have Hausdorff dimension zero. In such cases, a finer scale of measurement is required. To achieve this, one can introduce the generalized Hausdorff measure

\[
\mathcal{H}^h(E) = \lim_{\delta \to 0} \inf \left\{ \sum_i h(\operatorname{diam}(U_i)) : E \subset \bigcup_i U_i,\, \operatorname{diam}(U_i) < \delta \right\},
\]

where $h:(0,\infty)\to(0,\infty)$ is a nondecreasing and right continuous function such that $h(r)\to 0$ as $r\to 0$. Such a function is called a \emph{dimension function} or \emph{gauge function}. See \cite{Rogers} for a comprehensive treatment of these (outer) measures.

For the special choice $h(r)=r^s$, one recovers the classical $s$-dimensional Hausdorff measure. By allowing more general gauge functions, one can distinguish between sets of Hausdorff dimension zero that still exhibit substantial size or complexity when measured at a finer scale.

This approach has proved particularly useful for exceptional sets arising in Diophantine approximation and dynamical systems (for example \cite{MassTransferencePrinciple}), by allowing one to determine for which gauge functions $h$ the Hausdorff $h$--measure $\mathcal{H}^h$ of a set of zero Hausdorff dimension is positive, infinite, or zero. A classical example is the set of Liouville numbers $\L$. Recall that a real number $x \in \mathbb{R}\setminus \mathbb{Q}$ is called a \emph{Liouville number} if, for every $n \in \mathbb{N}$ there exist integers $p$ and $q$ with $q>1$ such that
\[
\left| x - \frac{p}{q} \right| < \frac{1}{q^n}.
\]

Olsen and Renfro determined the exact cut point between gauge functions $h$ for which $\mathcal{H}^h(\L)$ is infinite and those for which it is zero (see \cite{Ol1} and \cite{Ol2}). Furthermore, they proved that if $\mathcal{H}^h(\L)=\infty$ for a given gauge function $h$, then $\L$ does not have $\sigma$-finite $\mathcal{H}^h$ measure.

At the other extreme, for badly approximable numbers, another important class of sets arising in Diophantine approximation consists of numbers sharing the same Diophantine approximation constant. That is, for $x \in \mathbb{R} \setminus \mathbb{Q}$, the \emph{best constant of Diophantine approximation} of $x$ is

    $$
    k(x)=\sup\left\{ c>0 \mid \left| x - \frac{p}{q} \right| < \frac{1}{cq^2} \text{ has infinitely many solutions } \frac{p}{q}\in\Q\right\}
    $$

where $k$ can be regarded as a map from $\R\backslash\Q$ to $(0,\infty]$. The image $\{k(x):x\in\R\setminus\Q,k(x)<\infty\}$ of this function is known as the \emph{Lagrange spectrum}.

It is well known that $k^{-1}(3)$ is an uncountable dense set and it is a consequence of \cite[Theorem 1]{M:geometric_properties_Markov_Lagrange} that $k^{-1}(3)$ has Hausdorff dimension 0. Indeed, for $t\in\R$ define
\begin{multline}\label{eq:definition_K_t}
    K_t=\{[0;a_1,\dots,a_n,\dots]\ \mid\ \text{there exists $(a_{-n})_{n\geq 0}\in(\N_{>0})^{\N}$ such that} \\ 
[a_k;a_{k+1},\dots,]+[0;a_{k-1},a_{k-2},\dots]\leq t,\forall k\in\Z\}.
\end{multline}
Then $K_t\subset k^{-1}(-\infty,t]$ and $D(t)=\dim_H(K_t)=\overline{\dim_B}(K_t)=\dim_H(k^{-1}(-\infty,t])$ is a continuous function with $\max\{t\in\R:D(t)=0\}=3$.

In this paper, we determine the exact cut point at which the Hausdorff $h$--measure of $k^{-1}(3)$ drops from infinity to zero. This provides a complete characterization of all Hausdorff measures $\mathcal{H}^h(k^{-1}(3))$ without imposing any specific regularity assumptions on the dimension function $h$. Moreover, we show that if $\mathcal{H}^h(k^{-1}(3))=\infty$, then $k^{-1}(3)$ does not possess $\sigma$-finite $\mathcal{H}^h$ measure. Our main theorem is the following:

\begin{theorem}\label{mainth}
Let $h$ be a dimension function. 
\begin{enumerate}
\item If
\begin{equation*}
    \limsup_{\varepsilon\to 0}\frac{\log h(\varepsilon)}{\log\varepsilon}=0
\end{equation*}
then the set $k^{-1}(3)$ does not have $\sigma$-finite $\cH^h$ measure.
\item If
\begin{equation*}
    \limsup_{\varepsilon\to 0}\frac{\log h(\varepsilon)}{\log\varepsilon}>0
\end{equation*}
then $\cH^h(k^{-1}(3))=0$.
\end{enumerate}
Moreover, if $A_3$ consists of those $x\in k^{-1}(3)$ which are attainable, then $\cH^h(A_3)=\cH^h(k^{-1}(3))$ for any dimension function $h$.
\end{theorem}

\begin{definition}
An irrational number $x\in k^{-1}(3)$ is called \emph{attainable} if $ \left| x-\frac{p}{q} \right| \leq \frac{1}{3q^2} $ has infinitely many rational solutions $\frac{p}{q}$. 
\end{definition}

In comparison, one can say that the set of Liouville numbers $\L$ is strictly larger than $k^{-1}(3)$, in the sense that any dimension function $h$ for which $\mathcal{H}^h(k^{-1}(3)) = \infty$ also satisfies $\mathcal{H}^h(\L) = \infty$, whereas the converse does not necessarily hold. We will provide an explicit example of a dimension function $h$ such that $\mathcal{H}^h(\L) = \infty$ while $\mathcal{H}^h(k^{-1}(3)) = 0$.

A more general property of Liouville numbers is that they are \textit{immeasurable}, that is, any translation invariant Borel measure on $\R$ either assigns zero measure to $\L$ or is not $\sigma$--finite on $\L$ \cite[Theorem 1.1]{NoBorelMeasure} (note that this includes the Hausdorff $h$--measures). The main result of \cite{NoBorelMeasure} applies to a broad class of sets, however, it does not apply to $k^{-1}(3)$ because of several reasons. Recall that set of Liouville numbers $\L$ forms a $G_\delta$ set that is invariant by rational translations (it is not difficult to see that $k^{-1}(3)$ is not invariant by such translations, see for instance \cite{Davenport} and \cite{MendesImprovement}). A big difference is that from the topological point of view, the Liouville numbers are also larger than $k^{-1}(3)$. 

\begin{theorem}
The set $k^{-1}(3)$ is not $G_\delta$ and not $F_\sigma$, but it is $F_{\sigma\delta}$.
\end{theorem}

On the other hand, there are certain subsets of $k^{-1}(3)$ that are natural to consider. The use of generalized Hausdorff measure is then necessary to distinguish the size of these subsets. For example, the set $K_3$ is a Cantor set with Hausdorff dimension 0, since $K_3\subset k^{-1}(-\infty,3]$. For $t>3$, the Cantor set $K_t$ has the same positive Hausdorff dimension as $k^{-1}(-\infty,t]$, so it is natural to ask whether $K_3$ is as large as $k^{-1}(-\infty,3]$ from the point of view of Hausdorff dimension. Since $k^{-1}(-\infty,3)$ is countable because of Markov's theorem \cite[Theorem 16]{Bombieri}, it suffices to compare $K_3$ with $k^{-1}(3)$. The following result shows that $K_3$ is, in fact, significantly smaller than $k^{-1}(3)$. 

\begin{theorem}\label{thm:dimension_K_3}
Let $h$ be a dimension function. 
\begin{enumerate}
\item If 
\begin{equation*}
    \limsup_{\varepsilon\to 0}\frac{|\log h(\varepsilon)|}{\log|\log\varepsilon|}<3,
\end{equation*}
then $\cH^h(K_3)=\infty$. 

\item If 
\begin{equation*}
    \limsup_{\varepsilon\to 0}\frac{|\log h(\varepsilon)|}{\log|\log\varepsilon|}>3,
\end{equation*}
then $\cH^h(K_3)=0$. 

\end{enumerate}
\end{theorem}

\Cref{thm:dimension_K_3} does not decide the critical case $\limsup_{\varepsilon\to0}\frac{|\log h(\epsilon)|}{\log|\log \epsilon|} = 3$. The following result is sharper at this critical exponent. 

\begin{theorem}\label{thm:more_precise}
For a large enough constant $\hat{C}>0$, the dimension function 
\begin{equation*}
    g(\varepsilon)=\frac{\exp(\hat{C}\sqrt{\log|\log\varepsilon|})}{|\log\varepsilon|^3}
\end{equation*}
satisfies $\cH^g(K_3)=\infty$.

For the dimension function 
\begin{equation*}
    g(\varepsilon)=\frac{\log|\log\varepsilon|}{|\log\varepsilon|^3(\log\log|\log\varepsilon|)},
\end{equation*}
we have that $\cH^{g}(K_3)<\infty$.
\end{theorem}

While the first assertion of \Cref{thm:dimension_K_3} follows from the first assertion of \Cref{thm:more_precise}, we remark that there are functions $h$ that satisfy the second assertion of \Cref{thm:dimension_K_3} but are such that $h(\varepsilon)>\frac{\log|\log\varepsilon|}{|\log\varepsilon|^3(\log\log|\log\varepsilon|)}$ infinitely many times as $\varepsilon\to 0$. For example, this happens for the function
\begin{equation*}
    h(\varepsilon)=\exp\left(-x\left(3+\sin(\log x)\right)\right), \quad x = \log|\log\varepsilon|.
\end{equation*}

A consequence of the previous theorem, for example, is that for the dimension function $h(\varepsilon)=\frac{1}{|\log\varepsilon|^3}$ one has that $\cH^{h}(k^{-1}(3))=\infty$ but $\cH^{h}(K_3)=0$. In conclusion, there is a large gap between the generalized Hausdorff measure of these sets, so the set $K_3$ is significantly smaller than $k^{-1}(3)$.

Denote by $B_3$ the set of non--attainable elements of $k^{-1}(3)$, that is
\begin{equation*}
    B_3 = \left\{x\in\R\setminus\Q:\left|x-\frac{p}{q}\right|<\frac{1}{3q^2} \text{  has only finitely many solutions }\frac{p}{q}\in\Q\right\}.
\end{equation*}

A recent work \cite{cao2026badlyapproximablenumbers} by the first two authors with Zhe Cao imply the following result: for any $x=[x_0;x_1,x_2,\dots]\in B_3$, there is an $N\in\N$ such that $[0;x_{N+1},x_{N+2},\dots]\in K_3$. However, while the Cantor set $K_3$ is forward invariant by the Gauss map, the set $B_3$ is not. The next theorem shows that the non--attainable elements $B_3$ are even smaller than $K_3$ and thus smaller than $k^{-1}(3)$. 

\begin{theorem}\label{thm:very_badly_approximable_case}
Let $h$ be a dimension function. 
\begin{enumerate}
\item If 
\begin{equation*}
    \limsup_{\varepsilon\to 0}\frac{|\log h(\varepsilon)|}{\log|\log\varepsilon|}<2,
\end{equation*}
then $\cH^h(B_3)=\infty$. 

\item If 
\begin{equation*}
    \limsup_{\varepsilon\to 0}\frac{|\log h(\varepsilon)|}{\log|\log\varepsilon|}>2,
\end{equation*}
then $\cH^h(B_3)=0$. 
\end{enumerate}
\end{theorem}

The proof of the previous theorem is essentially the same as \Cref{thm:dimension_K_3}. The key distinction is that, in the case of $K_3$, one must account for permutations of renormalizable words, reflecting its invariance under the Gauss map. By contrast, the structure of $B_3$ is considerably more rigid so this additional counting is not required. For completeness, one could also formulate an analogue of \Cref{thm:more_precise} for $B_3$, but we refrain from doing so here in order to avoid unnecessary repetition. From this version it follows that $\cH^{h}(B_3)=0$ for the gauge function $h(\varepsilon)=\frac{1}{|\log\varepsilon|^2}$.

\section{Preliminaries}\label{sec:preliminaries}

\subsection{Continued fraction and Hausdorff dimension}

 Given an nondecreasing and right continuous function $h:[0,\infty)\to[0,\infty)$ with $h(0)=0$, we define the corresponding Hausdorff measure $\cH^h$ as
\begin{equation*}
    \cH^h(E) := \lim_{\delta\to0}\inf\left\{\sum h(\diam E_i): E_i \text{ open}, E\subset \bigcup E_i, \diam E_i \leq\delta\right\}.
\end{equation*}

We say that a set $E$ is an $h$--set if $0<\cH^h(E)<\infty$. It is not always possible to find a dimension function for which this holds.

Given $\underline{a}\in(\N_{>0})^\Z$ define
\begin{equation*}
    \lambda_i((a_n)_{n\in\Z})=[a_i;a_{i+1},a_{i+2},\dots]+[0;a_{i-1},a_{i-2},\dots].
\end{equation*}
By abuse of notation, if $(a_n)_{n\geq 0}$ is a one-sided sequence of positive integers, we also denote
\begin{equation*}
    \lambda_i((a_n)_{n\geq 0})=[a_i;a_{i+1},a_{i+2},\dots]+[0;a_{i-1},a_{i-2},\dots,a_1].
\end{equation*}

Continued fractions are useful to compute $k(x)$. Indeed, it is a classical fact that if $x=[a_0;a_1,a_2,\dots]$ then
\begin{equation}\label{eq:basic_identity}
    \left|x-\frac{p_n}{q_n}\right|=\frac{1}{(\gamma_{n+1}+\eta_{n+1})q_n^2},
\end{equation}
where $\gamma_{n+1}=[a_{n+1};a_{n+2},\dots]$ and $\eta_{n+1}=[0;a_n,a_{n-1},\dots,a_1]$. In particular we have the formula  $k(x)=\limsup_{n\to\infty}\lambda_n(x)$.

A curious and crucial identity to study irrational numbers $x$ with values $k(x)$ near 3 has been the following:
\begin{equation}\label{eq:cont_frac_identity}
    [2;2,z]+[0;1,1,z]=3, \qquad \text{for all $z>0$}.
\end{equation}

Given a finite word of positive integers $w=(a_1,\dots,a_n)\in(\N_{>0})^n$, define the closed sub-interval of $[0,1]$
\begin{equation*}
    I(w) := \{x\in[0,1]\ \mid\ x=[0; a_1,a_2,\dots,a_n,y], y\ge 1\} \cup \{[0,a_1,a_2, \dotsc,a_n]\},
\end{equation*}
consisting of the numbers in $[0,1]$ whose continued fractions start with $w$.

The precise diameter of such intervals is given by
\begin{equation*}
    \diam I(w)=\frac{1}{q_n(q_n+q_{n-1})},
\end{equation*}
where $q_k$ is defined recursively by $q_0=1$, $q_1=a_1$ and $q_{k+2}=a_{k+2}q_{k+1}+q_k$.

In particular, if $a_1,\dots,a_n\in\{1,2\}$, is easy to see that
\begin{equation}\label{eq:length_of_intervals_12}
    (3+2\sqrt{2})^{-n-1}<\diam I(w) < \left(\frac{3+\sqrt{5}}{2}\right)^{-n+1}
\end{equation}

The well known estimate of Jarnik \cite{Jarnik1928} says that the Hausdorff dimension of the Gauss-Cantor set $C(n)=\{[0;a_1,a_2,\dots]:1\leq a_i\leq n,\forall i\geq 1\}$ satisfies $1-C_2/n<\dim_H(C(n))<1-C_1/n$ for some explicit constants $0<C_1<C_2$ and all $n$ large. For our purposes, Jarnik result is enough, however a much more precise asymptotic estimate is due to Hensley \cite{Hensley2}:
\begin{equation}\label{eq:Hensley}
    \dim_H(C(n))=1-\frac{6}{\pi^2n}-\frac{72\log n}{\pi^4n^2}+O\left(\frac{1}{n^2}\right).
\end{equation}

See also \cite{DFSU} for an even more precise asymptotic estimate.

Another related result by Hensley estimates the cardinality of the set
\begin{equation*}
    Q_n(x)=\left\{(a,b)\in\Z^2: \substack{\displaystyle 1\leq a\leq b\leq x, \gcd(a,b)=1, \text{and all digits in} \\ \displaystyle \text{the continued fraction expansion of $a/b$ are $\leq n$}}\right\}.
\end{equation*}

Hensley \cite{Hensley2} proved that for any $n>1$, there is a constant $E_n>0$ such that $\lim_{x\to\infty}x^{-2\dim_H(C(n))}|Q_n(x)|=E_n$.  For our purposes, we will need a more elementary estimate that is not dependent on the rate of convergence.
\begin{proposition}
For all $x>1$,
\begin{equation}\label{eq:Q_n(x)}
    x^{2\dim_H(C(n))}\leq Q_n(x)\leq2x^{2\dim_H(C(n))}.
\end{equation}
\end{proposition}

\begin{proof}
Recall that a Gauss-Cantor set $K(B)$ is defined by an alphabet $B=\{\beta_1,\dots,\beta_\ell\}$ and for each word $\beta_j\in(\N_{>0})^{r_j}$, the interval $I_j=I(\beta_j)$ and $\psi|_{I_j}:= G^{r_j}|_{I_j}$ is an iterate of the Gauss map $G(x)=\frac{1}{x}-\lfloor \frac{1}{x}\rfloor$. This defines an expanding map $\psi\colon I_1\cup\dotsb\cup I_\ell\to I$ with Markov partition $I_1,\dots,I_\ell$. According to Palis--Takens \cite[Chapter 4, Pages 68-71]{PalisTakens}, if we let 
\[
    \lambda_j=\inf |\psi'|_{I_j}|, \qquad \Lambda_j=\sup |\psi'|_{I_j}|
\]
and $d_1,d_2\geq 0$ be such that
\[
    \sum_{i=1}^{\ell}\lambda_j^{-d_2}=1, \qquad \sum_{i=1}^{\ell}\Lambda_j^{-d_1}=1,
\]
then 
\begin{equation}\label{eq:Palis-Takens}
    d_1\leq \dim_H(K(B))\leq d_2.
\end{equation}

Let us discuss how to find estimates for $d_1$ and $d_2$. The iterates of the Gauss map are given explicitly by
\[
    \psi|_{I_j}(x)=\dfrac{q^{(j)}_{r_{j}}x - p^{(j)}_{r_{j}}}{p^{(j)}_{r_j-1} - q^{(j)}_{r_j-1}x}
\]
 where $\dfrac{p^{(j)}_k}{q^{(j)}_k} = [0; b^{(j)}_1, \dots, b^{(j)}_k]$ and $\beta_j = (b^{(j)}_1, \dots, b^{(j)}_{r_j})$. 

\noindent Hence 
\[
    (\psi|_{I_j})'(x)=\frac{(-1)^{r_j-1}}{(p^{(j)}_{r_j-1} - q^{(j)}_{r_j-1}x)^2}.   
\]

Let $x=[c_0,c_1,c_2,\dots]$ and $\frac{p_n}{q_n}=[c_0,c_1,\dots, c_n]$. Then
\[
    \frac{1}{2q_nq_{n+1}}<\frac{1}{q_n(q_n+q_{n+1})}<\left|x-\frac{p_{n}}{q_{n}}\right|<\frac{1}{q_nq_{n+1}},
\]
and therefore 
\[
    \frac{1}{2q_{n+1}}<|q_n x-p_n|<\frac{1}{q_{n+1}}.
\]

Therefore, 
\[
    (q_{r_j}^{(j)})^2<|(\psi|_{I_j})'(x)|=\frac{1}{(p^{(j)}_{r_j-1} - q^{(j)}_{r_j-1}x)^2}<(2q_{r_j}^{(j)})^2.
\]
Thus
\[
    (q_{r_j}^{(j)})^2\leq\lambda_j=\inf|\psi'|_{I_j}|\leq\Lambda_j=\sup|\psi'|_{I_j}|\leq (2q_{r_j}^{(j)})^2.
\]

Consider the set $P_n(x,m)=\{(a,b)\in Q_n(x):x/2^{m/2}<b\leq x/2^{(m-1)/2}\}$. Notice that writing $a/b=[0;b_1,\dots,b_\ell]$, each $[0;b_1,\dots,b_i]$ belongs to a different $P_n(x,m)$ for $1\leq b_i\leq\ell$. In particular the intervals $I(b_1,\dots,b_\ell)$ where $(a,b)\in P_n(x,m)$ form a Markov partition of the Gauss-Cantor set $C(n)$. By \eqref{eq:Palis-Takens} we have that 
\begin{equation*}
    1=\sum_{i=1}^{|P_n(x,m)|}\lambda_j^{-d_2}\leq |P_n(x,m)|x^{-2d_2}2^{d_1m}\leq |P_n(x,m)|x^{-2\dim_H(C(n))}2^{m},
\end{equation*}
whence
\begin{equation*}
    |Q_n(x)|=\sum_{m\geq 1}|P_n(x,m)|\geq x^{2\dim_H(C(n))}.
\end{equation*}
Similarly
\begin{equation*}
    1=\sum_{i=1}^{|P_n(x,m)|}\Lambda_j^{-d_1}\geq  |P_n(x,m)|x^{-2d_1}2^{d_2(m+1)}\geq  |P_n(x,m)|x^{-2\dim_H(C(n))}2^{m+1},
\end{equation*}
whence
whence
\begin{equation*}
    |Q_n(x)|=\sum_{m\geq 1}|P_n(x,m)|\leq 2x^{2\dim_H(C(n))}.
\end{equation*}
\end{proof}

\subsection{The subshift $\Sigma_3$}\label{subsec:Sigma_3}

Given $t>0$, define the subshift
\begin{equation*}
    \Sigma_t = \left\{\underline{a}\in (\N_{>0})^\Z \ \mid\ \sup_{n \in \Z} \lambda_n(\underline{a}) \leq t\right\}.
\end{equation*}

Note that for $t\in(0,\infty)$, one only needs finitely many symbols because $\Sigma_t\subset\{1,\dots,\lfloor t\rfloor\}^{\Z}$. Moreover, for $t=\sqrt{12}$ it is known that $\Sigma_t=\{1,2\}^\Z$ (see \cite[Chapter 1, Proposition 1.7]{LMMSbook}).

In general, the set $K_t$ is always closed, but it is not necessarily a Cantor set. For instance, when $t=\sqrt{13}$ the set $K_t$ includes the isolated point $[0;\overline{3}]$ (see \cite[Chapter 1, Lemma 7]{CF} and \cite[Chapter 1, Lemma 9]{CF}). More examples occur for all $\sqrt{5}\leq t<3$, where the set $K_t$ is countable (because of Markov's theorem). More generally, $K_t$ is a Cantor set whenever $\Sigma_t$ is a transitive subshift. Situations where the subshift $\Sigma_t$ is not transitive are for example when $t$ lies at the right endpoint of a gap of the Markov spectrum. Another particular situation when $\Sigma_t$ is not transitive is precisely at $t=3$: the periodic orbits $\dots 1111\dots$ and $\dots 2222\dots$ are both contained in $\Sigma_3$, yet there is no orbit in $\Sigma_3$ that contains both finite subwords $1111$ and $2222$ (see \cite[Lemma 11]{Bombieri}). However, the set $K_3$ is a Cantor set, as we will see.

The subshift $\Sigma_3$ was characterized in \cite{Reutenauer2006}: after substituting $22\mapsto a$ and $11\mapsto b$ it coincides with the subshift of all balanced bi-infinite words (eventually periodic and aperiodic). It is well known that the closure of any aperiodic Sturmian bi-infinite word gives a minimal infinite subshift (see \cite[Exercise 6.2.10]{Fogg}), so in particular, any orbit in that closure is an accumulation point. On the other hand, the eventually periodic balanced bi-infinite words can be approximated by the aperiodic (in the shift topology), thus any orbit of $\Sigma_3$ is an accumulation point (see \cite[Exercise 6.1.26]{Fogg} for example). Indeed, every finite factor of an eventually periodic balanced bi-infinite word occurs in an aperiodic Sturmian word, so the former can be approximated by the latter in the shift topology. In particular the shift $\Sigma_3$ is homeomorphic to a Cantor space and consequently the set $K_3$ is a Cantor subset of $\R$.

The shift $\Sigma_3$ can be described via the $S$--adic representation of balanced sequences. Background material and proofs of the statements below can be found in \cite{Reutenauerbook}. We call the substitutions $U$ and $V$ the \emph{inner renormalization operators}, defined by
\begin{align*}
    U : 
    \begin{aligned}
        &a\mapsto ab \\
        &b\mapsto b  \\
    \end{aligned}
    ,\quad
    V : 
    \begin{aligned}
        &a\mapsto a \\
        &b\mapsto ab.  \\
    \end{aligned}
\end{align*} 
Let $(\alpha,\beta)$ be a pair of finite words over the alphabet $\{a,b\}$. Define the \emph{exterior renormalization operators} $\overline{U},\overline{V}$ by
$$
\overline{U}: (\alpha,\beta)\mapsto (\alpha\beta,\beta),\quad \overline{V}: (\alpha,\beta)\mapsto (\alpha,\alpha\beta).
$$

The \emph{Cohn's tree} is a binary tree of finite words over $\{a,b\}$ defined as follows: take as root the word $ab$ and given a vertex, its children in the next layer are obtained by applying $U$ and $V$. Let $P$ be the set of vertices of Cohn's tree. These finite words together with $\{a,b\}$ are called \emph{(lower) Christoffel words}. It is a fact that a subword is a subfactor of a balanced sequence if and only if is a subfactor of some Christoffel word if and only if is a subfactor of some Sturmian sequence. By \cite[Lemma 3.11]{EGRS2024} we have a control of how large is this Christoffel word.

\begin{lemma}\label{lem:3w}
Let $w$ be a factor of a word in $P$. Then the length of the shortest word in $P$ containing $w$ is strictly smaller than $3|w|$.
\end{lemma}

Analogously, one can define a \emph{tree of alphabets} from word pairs as follows: the root is $(a,b)$ and given a vertex, its children in the next layer are obtained by applying $\overline{U}$ and $\overline{V}$. The set of vertices of this tree $\overline{P}$ are also called \emph{Christoffel pairs}.

The Cohn's tree and the tree of alphabets are related through the concatenation operation: given an alphabet $(\alpha,\beta)\in\overline{P}$ we have that $\alpha\beta\in P$ and conversely, given any Christoffel word $w\in P$, there is a unique factorization $w=\alpha\beta$ where $(\alpha,\beta)\in\overline{P}$. Moreover, given any $(\alpha,\beta)\in\overline{P}$, there are some substitutions $R_1,\dots,R_n$ such that $\alpha=(R_1\dots R_n)(a)$ and $\beta=(R_1\dots R_n)(b)$.

Another useful property of lower Christoffel words $P$ is that they are completely determined by counting their letters. More precisely, for $w\in P$, the \emph{Frobenius coordinates}  $(|w|_a+|w|_b,|w|_b)$ uniquely determine $w$, where $|w|_c$ denotes the number of occurrences of the letter $c\in\{a,b\}$. Furthermore, every word in $P$ is primitive, that is, no $w\in P$ is a proper power of a shorter word. Consequently, all cyclic permutations of a given $w$ are distinct. In addition, cyclic permutations of two distinct words $w, w' \in P$ cannot coincide, since their Frobenius coordinates are different.

Continued fractions are closely related to the Christoffel pairs. Indeed, for each continued fraction $[0;a_1,\dots,a_n]$ we can associate a unique Christoffel pair as follows. Define $(\alpha_0,\beta_0)=(a,b)$, $(\alpha_1,\beta_1)=\overline{V}^{a_1-1}(\alpha_0,\beta_0)$, $(\alpha_2,\beta_2)=\overline{U}^{a_2}(\alpha_1,\beta_1)$, $(\alpha_3,\beta_3)=\overline{V}^{a_3}(\alpha_2,\beta_2)$ and so on. In particular we will have that $q_0=|\alpha_0|/2$, $q_1=|\beta_1|/2$, $q_2=|\alpha_2|/2$, $q_3=|\beta_3|/2$ and so on, where $q_1,\dots,q_n$ are the denominators of the convergents of $[0;a_1,\dots,a_n]$. This process can be imagined as the mediant process used to build the convergents of the continued fraction $[0;a_1,\dots,a_n]$.

\begin{figure}[H]
\centering
\begin{tikzpicture}

\coordinate (L) at (0,0);
\coordinate (R) at (10,0);

\draw[-] (-0.5,0) -- (10.5,0);

\draw (L) -- ++(0,0.2);
\draw (R) -- ++(0,0.2);

\node[below] at (L) {$\alpha$};
\node[below] at (R) {$\beta$};


\foreach \k in {1,2,3,4,5}
{
    \coordinate (M\k) at ($ (L)!\k/6!(R) $);
    \draw (M\k) -- ++(0,0.2);
}

\node[below] at ($ (L)!1/6!(R) $) {$\alpha\beta$};
\node[below] at ($ (L)!2/6!(R) $) {$\alpha\beta^2$};
\node[below] at ($ (L)!3/6!(R) $) {$\alpha\beta^3$};
\node[below] at ($ (L)!4/6!(R) $) {$\alpha\beta^4$};
\node[below,yshift=-0.2cm] at ($ (L)!5/6!(R) $) {$\dots$};

\end{tikzpicture}

\begin{tikzpicture}

\coordinate (L) at (0,0);
\coordinate (R) at (10,0);

\draw[-] (-0.5,0) -- (10.5,0);

\draw (L) -- ++(0,0.2);
\draw (R) -- ++(0,0.2);

\draw (2,0.2) -- (2,0);
\draw (3,0.2) -- (3,0);
\draw (5,0.2) -- (5,0);
\draw (8,0.2) -- (8,0);
\draw (6.9,0.2) -- (6.9,0);

\node[below] at (L) {$\alpha$};
\node[below] at (R) {$\beta$};


\node[below] at (2,0) {$\alpha\beta^3$};
\node[below] at (3,-0.2) {$\dots$};
\node[below] at (5,0) {$\alpha\beta^3((\alpha\beta^3)^2\beta)^4$};
\node[below] at (6.9,-0.2) {$\dots$};
\node[below] at (8.3,0) {$(\alpha\beta^3)^2\beta$};

\end{tikzpicture}
\caption{Example of how constructing alphabets is similar to build the convergents of a number $[0;1,3,2,4]$.}
\end{figure}

Recall that given fractions $0<\frac{p}{q}<\frac{r}{s}$ the mediant is
\begin{equation*}
    \frac{p}{q}\oplus\frac{r}{s}=\frac{p+r}{q+s}.
\end{equation*}
Given $n\geq 1$, the Farey sequence $\mathcal{F}_n$ of order $n$ is the set of reduced coprime fractions $0\leq\frac{p}{q}\leq 1$ with $1\leq q\leq n$ increasingly ordered. Two fractions $\frac{p}{q}<\frac{r}{s}$ are \emph{Farey neighbors} in some $\mathcal{F}_n$ if and only if $ps-qr=1$ and $\max\{q,s\}\leq n<q+s$. The mediant $\frac{p+r}{q+s}$ will appear between $\frac{p}{q}$ and $\frac{r}{s}$ precisely at term $\mathcal{F}_{q+s}$.

The rule used to construct the tree of alphabets is related to the mediant rule used in the construction of the Farey sequence. We can build a bijection $\Theta$ from the Christoffel words $P\cup\{a,b\}$ to the positive fractions $\Q\cap[0,1]$ as follows. Given $w\in P\cup\{a,b\}$, we define $\Theta(w)=\frac{|w|_b}{|w|_{a,b}}$ where $|w|_{a,b}$ is the length of $w$ over the alphabet $\{a,b\}$ and $|w|_b$ is the number of letters $b$ in $w$. The following fact is consequence of \cite[Lemma 3.1]{GuguVillamil}.

\begin{lemma}\label{lem:farey_neighbors}
For each Farey neighbors $\frac{p}{q}<\frac{r}{s}$, there is a unique pair of words $(\alpha,\beta)\in\overline{P}$ such that $\Theta(\alpha)=p/q$ and $\Theta(\beta)=r/s$. 
\end{lemma}

We define $\Sigma(t, n)$ to be the set of length-$n$ subwords of sequences in $\Sigma_t$. From \cite[Theorem 1.1]{EGRS2024} we have the following:
	
\begin{theorem} \label{thm:equalities}
For all $n\geq 68$ we have
    \[
        \Sigma(3 + 6^{-3n}, n) = \Sigma(3, n) = \Sigma(3 - 6^{-3n}, n).
    \]
\end{theorem}

It is well known \cite{countingSturmianwords} that the cardinality of all Sturmian factors is asymptotically equivalent to $n^3/\pi^2$. In particular, we have a good control on the cardinality of the subfactors of $\Sigma_3$.

\begin{lemma}\label{lem:counting_sturminan_words}
$|\Sigma(3,n)|\sim n^3/(4\pi^2)$.
\end{lemma}
\begin{proof}
Let $L_n$ denote the set of all Sturmian factors of length $n$. From \cite{countingSturmianwords} we know that $|L_n|\sim n^3/\pi^2$. We will prove the following:

\begin{equation*}
|\Sigma(3,n)| =
\begin{cases}
|L_{n/2}| + |L_{n/2+1}| - 2, & n \text{ even},\\[2pt]
2\,|L_{(n+1)/2}| - 2, & n \text{ odd}.
\end{cases}
\end{equation*}

First suppose $n$ even. Given $w\in\Sigma(3,n)$, if $w$ begins with an even block of 1's or 2's then by applying the inverse substitution $22\mapsto a$, $11\mapsto b$ it can be written as an Sturmian factor of length exactly $n/2$, so we have $|L_{n/2}|$ such words. If it begins with an odd block, then it also ends with an odd block so by adding one digit $\{1,2\}$ (which are uniquely determined because $w$ contains both 1 and 2 in this case) at each side and applying the substitution, we obtain a non constant Sturmian word of length $n/2+1$, so we obtain $|L_{n/2+1}|-2$ factors.

Now assume $n$ odd. In this case the elements of $\Sigma(3,n)$ can be obtained from $L_{(n+1)/2}$ by substituting $a\mapsto 22$, $b\mapsto 11$ and then deleting precisely the first or the last digit. This gives two different factors for each $L_{(n+1)/2}$ except for the constant factors.
\end{proof}

\begin{remark}
In fact, it is possible to prove that $|\Sigma(3,n)|\leq 9n^3$ for all $n\geq 1$, see \cite[Corollary 3.13]{EGRS2024}.
\end{remark}

\begin{lemma}\label{lem:existence_of_n0}
Let $t\in L$ and $\varepsilon>0$. Given any $x=[0;a_1,a_2,a_3,\dots]\in k^{-1}(t)\cap(0,1)$, there is a $n_0$ and a sequence of positive integers $(b_j)_{j<n_0}$ such that 
\begin{equation*}
    (\dots,b_{n_0-2},b_{n_0-1},a_{n_0},a_{n_0+1},a_{n_0+2},\dots)\in\Sigma_{t+\varepsilon}.
\end{equation*}
\end{lemma}
\begin{proof}
First, note that for $n$ large $a_n,a_{n+1},\dots\in\{1,2,\dots,\lfloor t\rfloor\}$. Since $k(x)=\limsup_{i\to\infty}\lambda_i(x)=t$, there is a subsequence $(i_k)_k$ such that $\lim_{k\to\infty}\lambda_{i_k}(x)=t$. In particular, by taking a further subsequence, we can assume that the subwords $(a_{i_k}),(a_{i_{k+1}-1},a_{i_{k+1}},a_{i_{k+1}+1}), (a_{i_{k+2}-2},a_{i_{k+2}-1},a_{i_{k+2}},a_{i_{k+2}+1},a_{i_{k+2}+2}),\dots$ converge to a bi-infinite sequence $\underline{b}=(b_n)_{n\in\Z}$ with $\sup_{n\in\Z}\lambda_i(\underline{b})=t$. Take $m=\lceil \log_2 (1/\varepsilon)\rceil+1$ and $n_0:=i_k$ so large such that $\lambda_{i}(x)<t+\varepsilon/2$ for all $i\geq i_k$ and such that  $(a_{i_k-m},\dots,a_{i_k},\dots,a_{i_k+m})=(b_{i_k-m},\dots,b_{i_k},\dots,b_{i_k+m})$. In particular we have that $(\dots,b_{n_0-2},b_{n_0-1},a_{n_0},a_{n_0+1},a_{n_0+2},\dots)$ has Markov value at most $t+\varepsilon/2+2^{-m}<t+\varepsilon$.
\end{proof}

\section{Dimensionality of $k^{-1}(3)$}

Recall that $k^{-1}(3)$ denotes the set of real numbers $x \in \mathbb{R}$ such that,
\[
\left| x - \frac{p}{q} \right| \leq \frac{1}{(3+\varepsilon)q^2}
\]
has infinitely many rational solutions $\frac{p}{q}$ for any $\varepsilon<0$ but only finitely many for any $\varepsilon>0$. An interesting subset of these numbers is the set of \emph{attainable} numbers, defined as follows.

\begin{definition}\label{def:attainable}
An irrational number $x\in k^{-1}(3)$ is called \emph{attainable} if $ \left| x-\frac{p}{q} \right| \leq \frac{1}{3q^2} $ has infinitely many rational solutions $\frac{p}{q}$. 
\end{definition}
Now we will show that there is a large class of numbers in $k^{-1}(3)$ which are attainable.
\begin{lemma}\label{lem:attainable_criteria}
Let $(e_i)_{i\geq 1}$ be any sequence of positive integers with $\lim_{i\to\infty}e_i=\infty$. Then the irrational number
\begin{equation*}
    x = [0;1^{e_1},2,2,1^{e_2},2,2,1^{e_3},2,2,\dots],
\end{equation*}
satisfies $k(x)=3$ and is attainable.
\end{lemma}

\begin{proof}
Let $i_0\geq 2$ be such that $e_i\geq 2$ for all $i\geq i_0$. If $e_{i+1}$ is even and $e_i$ odd, or if $e_{i}<e_{i+1}+2$ are both odd, then using \eqref{eq:basic_identity} one has
\begin{align}\label{eq:bad_cut1}
\begin{split}
    \lambda(\dots 221^{e_{i}}|221^{e_{i+1}}22\dots)&=[2;2,1^{e_{i+1}},2,2,\dots]+[0;1^{e_i},2,2,\dots] \\
    &= 3 + [0;1^{e_i},2,2,\dots] - [0;1^{e_{i+1}+2},2,2,\dots] > 3.
\end{split}
\end{align}
If $e_{i+1}$ is odd and $e_i$ even, or if $e_{i+1}>e_i+2$ are both even, or if $e_{i+1}<e_i+2$ are both odd, then using \eqref{eq:basic_identity} again
\begin{align}\label{eq:bad_cut2}
\begin{split}
    \lambda(\dots 221^{e_{i}}2|21^{e_{i+1}}22\dots)&=[0;1^{e_{i+1}},2,2,\dots]+[2;2,1^{e_i},2,2,\dots] \\
    &= 3 + [0;1^{e_{i+1}},2,2,\dots] - [0;1^{e_{i}+2},2,2,\dots] > 3.
\end{split}
\end{align}
Note that, disregarding parities, since $e_i\to\infty$ as $i\to\infty$, the equations \eqref{eq:bad_cut1} and \eqref{eq:bad_cut2} show that $k(x)=\limsup_{n\to\infty}\lambda_n(x)=3$.

Now we will show that for infinitely many positions $n$ one has $\lambda_n(x)>3$. First, assume that there infinitely many indices $i\geq i_0$ such that $e_i$ and $e_{i+1}$ have different parities. From both equations above we see that we have infinitely many positions $n$ such that $\lambda_n(x)=\gamma_{n+1}+\eta_{n+1}>3$.

Now, assume that after some larger $i_0$, all $e_i$ have the same parity for $i\geq i_0$. If such parity is odd, then using \eqref{eq:bad_cut1} and \eqref{eq:bad_cut2} it suffices to consider those infinitely many indices such that $e_i<e_{i+1}+2$ or $e_{i+1}<e_i+2$. If the parity is eventually even and $e_{i+1}>e_i+2$ infinitely many times, then using \eqref{eq:bad_cut2} we are done. Otherwise, the pattern $221^{2k}221^{2k+2}221^{2k+4}221^{2k+6}$ will repeat infinitely many times (with increasing $k$) and in this case we use \eqref{eq:basic_identity} again to conclude that
\begin{align*}
    &\lambda(\dots 221^{2k}221^{2k+2}2|21^{2k+4}221^{2k+6}\dots)\\
    &=[0;1^{2k+4},2,2,1^{2k+6},\dots]+[2;2,1^{2k+2},2,2,1^{2k},2,2,\dots] \\
    &=3+[0;1^{2k+4},2,2,1^{2k+6},\dots]-[0;1^{2k+4},2,2,1^{2k},2,2,\dots] > 3,
\end{align*}
infinitely many times.

\end{proof}

Dimension functions $h$ for which $\mathcal{H}^h(h^{-1}(3)) = \infty$ 
will be obtained by constructing suitable Cantor subsets 
$C \subset k^{-1}(3)$. 
In view of this, we introduce the following lemma concerning the fractal geometry of Cantor sets. 
The next lemma is a slightly modified version of \cite[Lemma~1]{Mauduit}. 
The proof is similar, and we include it here for completeness.

\begin{lemma}\label{cantor}
Let $K\subset\R$ be a Cantor set obtained as $K=\bigcap_{n\geq 1}K_n$ where $K_n$ is a union of intervals of size at least $\varepsilon_n>0$ with disjoint interiors such that, for each interval $I$ of $K_n$, we have 
$$|\left\{  J\text{ interval of } K_{n+1}|\ J\subset I \right\}|\geq m(n)>0.$$

Let $N_k=\Pi_{n<k}m(n)$, and let $h:[0,\infty)\to[0,\infty)$ be a dimension function such that $\lim_{k\to\infty}N_k\cdot h(\varepsilon_{k+1})=\infty$, then $\cH^h(K)=\infty$.
\end{lemma}

  \begin{proof}
Let us consider a finite covering of $K$ by intervals $K \subset \bigcup_{j=1}^{\ell} I_j$, and denote, for any $j \in \{1, \ldots, \ell\}$, by $v(j)$ the unique integer such that $\varepsilon_{v(j)+1} \leq |I_j| < \varepsilon_{v(j)}$.

We now replace each interval $I_j$ by at most two intervals of $K_{v(j)}$ whose union contains $I_j \cap K$. 
Note that this is possible because if $I_j$ intersects three or more intervals of $K_{v(j)}$, 
then $I_j$ must contain one of them, and hence $\varepsilon_{v(j)} < \operatorname{diam}(I_j)$.
 So, we get a covering of $K$ by intervals of size at least $\varepsilon_{v(j)}$ of $K_{v(j)}$, from which we can extract a sub-covering 
\[
K \subset \bigcup_{s=1}^{\tilde{\ell}} \tilde{I}_s
\]
by disjoint intervals $\tilde{I}_s$ of sizes at least $\varepsilon_{n(s)}$ of $K_{n(s)}$, 
and since each $I_j$ was subdivided into at most two intervals, 
we have, for every positive integer $n$,
\begin{equation}\label{eq4}
|\{1\leq s\leq\tilde{\ell} \mid n(s)=n\}| \le 2|\{1\leq j\leq\ell \mid v(j)=n\}|.
\end{equation}

\begin{claim}
If we have a finite collection of disjoint intervals $\tilde{I}_s$ covering $K$ such that, for all $s$, $\tilde{I}_s$ is an interval of $K_{n(s)}$ for some $n(s)$, and if we denote $k_n := |\{1\leq s\leq\tilde{\ell} \mid n(s)=n\}|$, then we have 
\[
\sum_{n\ge1} k_n / N_n \ge 1.
\]
\end{claim}

To prove the claim, let us first remark that it is trivial when $n(s)$ is constant and equal to $\tilde{n}$, since in this case the covering $\bigcup_{s=1}^{\tilde{\ell}}\tilde{I}_s$ coincides with $K_{\tilde{n}}$ which in turn is the union of at least $N_{\tilde{n}}$ disjoint intervals of size $\varepsilon_{\tilde{n}}$, so $\sum_{n\ge1} k_n / N_n = k_{\tilde{n}} / N_{\tilde{n}} \ge 1$.

Otherwise, let us take $n^*$ maximal, such that $k_{n^*} \ne 0$, and change all intervals of $K_{n^*}$ which belong to the collection and intersect a given interval $I$ of $K_{n^*-1}$ by $I$. This is possible because, whenever one interval of $K_{n^*}$ in the collection intersects $I$, all intervals of $K_{n^*}$ intersecting $I$ necessarily belong to the collection. By doing so, we replace at least $N_{n^*}/N_{n^*-1}$ intervals of $K_{n^*}$ by one interval of $K_{n^*-1}$, so that the positive integer $\sum_{n\ge1} n \cdot k_n$ diminishes and, since
\[
(k_{n^*} - N_{n^*}/N_{n^*-1})/N_{n^*} + (k_{n^*-1} + 1)/N_{n^*-1}
= k_{n^*}/N_{n^*} + k_{n^*-1}/N_{n^*-1},
\]
the value of $\sum_{n\ge1} k_n/N_n$ does not increase. We repeat this process until all intervals belong to the same $K_r$, when we have at least $N_r$ intervals. Since $N_r/N_r = 1$, the claim is proved.

Now we can conclude the proof of \Cref{cantor}.

As we have seen above, (some of) the intervals of size at least $\varepsilon_{v(j)}$ we considered form a covering of $K$ by disjoint intervals $\tilde{I}_s$ where, for all positive integers $1\leq s\leq\tilde{\ell}$, $\tilde{I}_s$ is an interval of $K_{n(s)}$ for some $n(s)$. It follows from the claim and from \eqref{eq4} that 
\[
2 \sum_{r \ge 1} (m_r / N_r) \ge 1,
\]
where we put $m_r = |\{ j \mid m(j) = r\}|$.

Finally, we have
\[
\sum_{j=1}^{\ell} h(|I_j|) 
\ge \sum_{j=1}^{\ell} h(\varepsilon_{v(j)+1}) 
\ge \sum_{r \ge r_0} v_r \cdot h(\varepsilon_{r+1}),
\]
where $r_0 = \min\{r \mid k_r > 0\}$.
Since $\sum_{r \ge r_0} (m_r / N_r) \ge 1/2$, it follows that
\[
\sum_{j=1}^{\ell} h(|I_j|) 
\ge \sum_{r \ge r_0} \frac{m_r}{N_r} \cdot N_r \cdot h(\varepsilon_{r+1})
\ge \frac{1}{2} \cdot \min_{r \ge r_0} (N_r \cdot h(\varepsilon_{r+1})),
\]
which tends to $+\infty$ when $r_0 \to +\infty$. But when the norm of the covering tends to $0$, we do have $r_0 \to +\infty$, and the proof is complete.
\end{proof}

We now proceed to prove \Cref{mainth}. The proof will be divided into the following three lemmas.

\begin{lemma}\label{lem:infinite_case}
Let $h$ be a dimension function. If
\[
\lim_{\varepsilon\to0}\frac{\log h(\varepsilon)}{\log\varepsilon}=0
\]
then $\mathcal H^h(k^{-1}(3))=\infty$.
\end{lemma}

\begin{proof}
The goal is to construct a Cantor set $C\subset k^{-1}(3)$ with $\mathcal H^h(C)=\infty$.
Let $c:=\log(3+2\sqrt2)>1$ and, for $T>0$, put
\[
\delta(T):=\sup_{0<\varepsilon\le e^{-T}}\frac{|\log h(\varepsilon)|}{|\log\varepsilon|}.
\]
Then $\delta$ is non-increasing, and $\delta(T)\to0$ as $T\to\infty$ by hypothesis.

\medskip

Fix $T_0>0$ such that $\delta(T_0)\le\frac1{20c}$, and let $f_1:=5$. Having defined $f_1,\dots,f_k$, put
\[
S_k:=f_1+\dots+f_k,\qquad T_k:=T_0+2c\,S_k,
\]
and
\[
f_{k+1}:=\begin{cases}
2f_k & \text{if } 8c\,f_k\,\delta(T_k+4cf_k)\le1,\\[2pt]
f_k & \text{otherwise.}
\end{cases}
\]
Thus $(f_k)$ is a non-decreasing sequence of integers $\ge5$. We claim that
\begin{equation}\label{eq:invariant}
4c\,f_k\,\delta(T_k)\le1\qquad\text{for all }k\ge1.
\end{equation}
For $k=1$ we have $T_1\ge T_0$, hence $4cf_1\delta(T_1)\le20c\,\delta(T_0)\le1$.
If $f_{k+1}=2f_k$, then $T_{k+1}=T_k+4cf_k$ and $4cf_{k+1}\delta(T_{k+1})=8cf_k\,\delta(T_k+4cf_k)\le1$ by the definition of $f_{k+1}$.
If $f_{k+1}=f_k$, then $T_{k+1}>T_k$, so $\delta(T_{k+1})\le\delta(T_k)$, and \eqref{eq:invariant} for $k$ implies it for $k+1$.

\medskip

Now we construct a Cantor set determined by the sequences defined above. For each \( k \), let \( r_i \in \mathbb{N} \) for \( i = 1, \dots, k \) satisfy
\[
f_i-4 \le r_i < 2f_i-4.
\]
We consider the set of real numbers whose continued fraction expansion begins with the following pattern:
\[
\underbrace{1\,1\,\cdots\,1}_{r_1\text{ times}}\,2\,2\,\underbrace{1\,1\,\cdots\,1}_{r_2\text{ times}}\,2\,2\,\cdots\,\underbrace{1\,1\,\cdots\,1}_{r_k\text{ times}}\,2\,2.
\]
That is, the expansion consists of alternating blocks of \( r_i \) ones followed by two 2's, for \( i = 1 \) to \( k \). This defines a finite-level construction of a Cantor set based on the growth of the sequences \( (f_i) \).
Let \( C \) be the Cantor set obtained by repeating the above construction for all \( n \in \mathbb{N} \). Then the set \( C \) satisfies the following properties:

\begin{enumerate}
    \item \( C \subset k^{-1}(3) \). This follows from the fact that $r_i\to\infty$. Therefore, any $x\in C$ is attainable and $k(x)=3$ because of Lemma \ref{lem:attainable_criteria}.
    \item At step \( n \) of the construction, the set \( C \) consists of exactly \( N_n:=\prod_{i=1}^{n} f_i \) closed intervals. This is because for each index \( i \), the repetition number \( r_i \) ranges over exactly \( f_i \) possible values, independent from the other levels.
    \item Each closed interval appearing in step \( n \) of the construction has length at least \( e^{-(2cf_1+\dots+2cf_n)}=e^{-2cS_n}\geq e^{-T_0-2cS_n}=: \varepsilon_n \). This is a classical fact in the theory of continued fractions (see \eqref{eq:length_of_intervals_12}).
\end{enumerate}

\medskip

Now, we will apply Lemma \ref{cantor} with $m(n)=f_{n+1}$. Since $f_1=5$ is a constant, the hypothesis $N_k\,h(\varepsilon_{k+1})\to\infty$ of that lemma is equivalent to
\[
f_1\cdots f_n\;h(\varepsilon_{n+1})\longrightarrow\infty .
\]
By the definition of $\delta$, $\log h(\varepsilon_{n+1})\ge-\delta(T_{n+1})\,T_{n+1}$, where
$T_{n+1}=T_0+2c\sum_{i=1}^{n+1}f_i$. For $i\le n+1$ we have $T_{n+1}\ge T_i$, hence by \eqref{eq:invariant}
\[
2c\,f_i\,\delta(T_{n+1})\le2c\,f_i\,\delta(T_i)\le\tfrac12 ,
\]
and also $\delta(T_{n+1})\,T_0\le\delta(T_0)\,T_0\le T_0/(20c)$. Therefore
\[
\log\bigl(f_1\cdots f_n\,h(\varepsilon_{n+1})\bigr)\;\ge\;\sum_{i=1}^n\log f_i-\frac{n+1}{2}-\frac{T_0}{20c}
\;\ge\;n\Bigl(\log5-\frac12\Bigr)-\frac12-\frac{T_0}{20c}\;\longrightarrow\;\infty ,
\]
because $\log5>\frac12$. By Lemma \ref{cantor}, $\mathcal H^h(C)=\infty$, and since $C\subset k^{-1}(3)$ we conclude that
$\mathcal H^h(k^{-1}(3))\ge\mathcal H^h(C)=\infty$.
\end{proof}

\begin{lemma}
If $\lim_{\varepsilon\to0}\frac{\log h(\varepsilon)}{\log \varepsilon}=0$, then $k^{-1}(3)$ does not have $\sigma$-finite $\cH^h$-measure.
\end{lemma}

\begin{proof}
Define the dimension function $g(x)=h(x)^2$. We have that
\begin{equation*}
    \lim_{\varepsilon\to0}\frac{\log g(\varepsilon)}{\log\varepsilon}=\lim_{\varepsilon\to0}\frac{2\log h(\varepsilon)}{\log\varepsilon}=0.
\end{equation*}
Hence $\cH^g(k^{-1}(3))=\infty$. Moreover
\begin{equation*}
    \lim_{\varepsilon\to0}\frac{g(\varepsilon)}{h(\varepsilon)}=\lim_{\varepsilon\to0}h(\varepsilon)=0.
\end{equation*}
In particular, if $k^{-1}(3)=\cup_n E_n$ with $\cH^h(E_n)<\infty$, we must have $\cH^g(E_n)=0$. But this would imply that
\begin{equation*}
    \cH^g(k^{-1}(3))=\cH^g(\cup_n E_n)\leq \sum_n\cH^g(E_n)=0,
\end{equation*}
which is clearly a contradiction.
\end{proof}

\begin{lemma}
Let $h$ be a dimension function. If 
\begin{equation*}
    \limsup_{\varepsilon\to 0}\frac{\log h(\varepsilon)}{\log\varepsilon}>0
\end{equation*}
then $\cH^h(k^{-1}(3))=0$.
\end{lemma}

\begin{proof}
By hypothesis, there is a positive $d$ and a decreasing sequence $\varepsilon_k$ such that $\varepsilon_k\to0$ when $k\to\infty$ such that
\begin{equation}\label{eq:fast_decay_h}
    h(\varepsilon_k)\leq \varepsilon_k^d, \quad \text{for all $k$}.
\end{equation}
Take $m\in\N$ so large so that $9m^3<2^{dm-1}$.

Define the sequence of positive integers $r_k=m\left\lceil\frac{|\log_2\varepsilon_k|}{m}\right\rceil$ for all $k\geq 1$.

By \Cref{lem:existence_of_n0}, we know there is $n_0=n_0(x)\geq 2$ and a sequence of positive integers $(b_j)_{j<n_0}$ such that 
\begin{equation}\label{eq:existence_of_n_0}
    (\dots,b_{n_0-2},b_{n_0-1},a_{n_0},a_{n_0+1},a_{n_0+2},\dots)\in\Sigma_{3+6^{-m}}.
\end{equation}

In particular, this shows that $k^{-1}(3)$ can be covered by the countable union
\begin{equation*}
    k^{-1}(3)\subset \bigcup_{n_0\geq 2}\bigcup_{(a_1,\dots,a_{n_0-1})\in(\N_{>0})^{n_0-1}}\{[0;a_1,\dots,a_{n_0-1},\gamma]:\gamma\in K_{3+6^{-m}}\}.
\end{equation*}
where $K_t$ was defined on \eqref{eq:definition_K_t}. If we show that $\{[0;a_1,\dots,a_{n_0-1},\gamma]:\gamma\in K_{3+6^{-m}}\}$ has zero Hausdorff $h$--measure, then $\cH^h(k^{-1}(3))=0$. From now on we fix $(a_1,\dots,a_{n_0-1})$ and assume that $x\in\{[0;a_1,\dots,a_{n_0-1},\gamma]:\gamma\in K_{3+6^{-m}}\}$.

It follows from \eqref{eq:existence_of_n_0} that given any sub-block of length $m$ at a position $n\geq n_0$ of $x=[0;a_1,a_2,\dots]$, that is, a word of the form $(a_n,\dots,a_{n+m-1})$, it must  belong to $\Sigma(3,m)$ and we have at most $9m^3<2^{md-1}$ possibilities for such sub-blocks.

We want to build a covering of $k^{-1}(3)$ with arbitrarily small $h$--measure. For this, we will employ \Cref{thm:equalities} to cover each $\{[0;a_1,\dots,a_{n_0-1},\gamma]:\gamma\in K_{3+6^{-m}}\}$.

Observe that we can divide the block $w=(a_{n_0},a_{n_0+1},\dots,a_{n_0+r_k-1})$ in the subblocks $(a_{n_0+im},\dots,a_{n_0+(i+1)m-1})$ for $i=0,\dots,r_k/m-1$. In particular, we can cover the block $B$ with at most 
\begin{align*}
    (2^{md-1})^{r_k/m}&=(2^{md-1})^{\lceil|\log_2\varepsilon_k|/m\rceil} \\
    &\leq (2^{md-1})^{|\log_2\varepsilon_k|/m+1} \\
    &=(2^{|\log_2\varepsilon_k|})^d\cdot 2^{-|\log_2\varepsilon_k|/m}\cdot 2^{md-1} \\
    &=\varepsilon_k^{-d}\cdot 2^{-|\log_2\varepsilon_k|/m}\cdot 2^{md-1}
\end{align*}
intervals $I(p,w)$ where $p=(a_1,\dots,a_{n_0-1})$ is fixed and $w\in\Sigma(3,n)^{r_k/m}$. The length of the intervals $I(p,w)$ is at most $2^{-n_0-r_k}\leq 2^{-r_k}\leq 2^{-|\log_2\varepsilon_k|}= \varepsilon_k$. Using \eqref{eq:fast_decay_h}, we have that the measure of the covering $I(p,w)$ is at most
\begin{align*}
    \sum_{w\in\Sigma(3,n)^{r_k/m}}h(\diam I(w))&\leq  \sum_{w\in\Sigma(3,n)^{r_k/m}}h(\varepsilon_k)\\
    &\leq (2^{md-1})^{r_k/m}\cdot\varepsilon_k^d \\
    &\leq\varepsilon_k^{-d}\cdot 2^{-|\log_2\varepsilon_k|/m}\cdot 2^{md-1}\cdot\varepsilon_k^d \\
    &=2^{md-1}\cdot2^{-|\log_2\varepsilon_k|/m}.
\end{align*}

Since the sequence $\varepsilon_k$ is decreasing and $2^{md-1}$ is constant, the above term goes to zero as $k\to\infty$, which shows that $\cH^h(k^{-1}(3))=0$.

\end{proof}
We would like to remark that our result can be compared with those obtained by Olsen and Renfro in \cite{Ol1, Ol2}, concerning the Hausdorff measure of the set of Liouville numbers $\L$. In these articles, the authors proved the following
\begin{theorem}\label{thm:liouville_measure}
Let $h$ be a dimension function. Then:
\begin{enumerate}
    \item If $\limsup_{r \to 0} \frac{\Gamma_h(r)}{r^t} = 0$ for some $t > 0$, then $\mathcal{H}^h(\L) = 0$.
    \item If $\limsup_{r \to 0} \frac{\Gamma_h(r)}{r^t} > 0$ for all $t > 0$, then the set $\L$ does not have $\sigma$-finite $\mathcal{H}^h$-measure.
\end{enumerate}
Here, $\Gamma_h(r) = r\cdot\inf_{0 < s \leq r} \frac{h(s)}{s}$ and $\L$ denotes the set of Liouville numbers.
\end{theorem}

It is proved in \cite[Lemma 2.2]{Ol2} that for all dimension functions $h$ and all sets $E \subset \mathbb{R}$, we have:
\[
\mathcal{H}^{\Gamma_h}(E) \leq \mathcal{H}^h(E) \leq 2 \mathcal{H}^{\Gamma_h}(E).
\]
Using this, we may reformulate \Cref{mainth} as follows,

\begin{theorem}
Let $h$ be a dimension function. Then:
\begin{enumerate}
    \item If $\liminf_{r \to 0} \frac{\Gamma_h(r)}{r^t} = 0$ for some $t > 0$, then $\mathcal{H}^h(k^{-1}(3)) = 0$.
    \item If $\liminf_{r \to 0} \frac{\Gamma_h(r)}{r^t} > 0$ for all $t > 0$, then the set $k^{-1}(3)$ does not have $\sigma$-finite $\mathcal{H}^h$-measure.
\end{enumerate}
\end{theorem}

This comparison shows the set of Liouville numbers is strictly larger than $k^{-1}(3)$ from the point of view of Hausdorff dimension, since there are dimension functions $h$ for which $\cH^{h}(\L)=\infty$ while $\cH^{h}(k^{-1}(3))=0$. Indeed, consider for example the right continuous dimension function
\begin{equation*}
    h(\varepsilon)=\begin{cases}
        e^{-(2k+1)!}\cdot\left(e^{(2k)!}\right)^{1-1/k}, & \text{if } e^{-(2k+2)!}\leq\varepsilon\leq e^{-(2k+1)!}; \\
        \varepsilon\cdot\left(e^{(2k)!}\right)^{1-1/k}, & \text{if } e^{-(2k+1)!}\leq\varepsilon< e^{-(2k)!}.
    \end{cases}
\end{equation*}

It is easy to see that
\begin{equation*}
    \Gamma_h(\varepsilon)=\begin{cases}
        \varepsilon\cdot\left(e^{(2k)!}\right)^{1-1/k}, & \text{if } e^{-(2k+1)!}\leq\varepsilon \leq \left(e^{(2k)!}\right)^{-1+1/k}\cdot h(e^{-(2k-1)!}); \\
        h(e^{-(2k-1)!}), & \text{if } \left(e^{(2k)!}\right)^{-1+1/k}\cdot h(e^{-(2k-1)!})\leq\varepsilon< e^{-(2k-1)!}.
    \end{cases}
\end{equation*}

In particular, for $\varepsilon_k=e^{-(2k+1)!}$ we have that 
\begin{equation*}
    \lim_{k\to\infty}\frac{\log h(\varepsilon_k)}{\log \varepsilon_k}=1-\lim_{k\to\infty}\frac{(2k)!(1-1/k)}{(2k+1)!}=1,
\end{equation*}
which implies $\limsup_{\varepsilon\to0}\frac{\log h(\varepsilon)}{\log\varepsilon}>0$ whence $\cH^{h}(k^{-1}(3))=0$ because of \Cref{mainth}. On the other hand, for $\varepsilon_k=e^{-(2k)!-1}$ and any $t>0$ one has
\begin{equation*}
    \lim_{k\to\infty}\frac{\Gamma_h(\varepsilon_k)}{\varepsilon_k^t}=\lim_{k\to\infty}e^{t-1}\left(e^{(2k)!}\right)^{1-1/k+t-1}=\lim_{k\to\infty}e^{t-1}\left(e^{(2k)!}\right)^{t-1/k}=\infty,
\end{equation*}
which implies $\limsup_{\varepsilon\to0}\frac{\Gamma_h(\varepsilon)}{\varepsilon^t}>0$ for all $t>0$, so $\cH^{h}(\mathbb{L})=\infty$ because of \Cref{thm:liouville_measure}.

Finally, let us mention that from the topological point of view, the Liouville numbers are also larger than $k^{-1}(3)$. First, we will show that $k^{-1}(3)$ is not $F_\sigma$. Indeed, given some knowledge about the generalized Hausdorff measures of a set implies some restriction on its topology, so the following lemma can be of independent interest.

\begin{lemma}\label{lem:general_topology}
Suppose that $E\subset\R$ is a bounded set of Hausdorff dimension $0$ such that for any dimension function $h$ with
\begin{equation}\label{eq:decay_ratio_general_topology}
    \lim_{\varepsilon\to0}\frac{\log h(\varepsilon)}{\log\varepsilon}=0,
\end{equation}
one has $\cH^{h}(E)>0$. Then $E$ is not $F_{\sigma}$.
\end{lemma}

\begin{proof}
Suppose by contradiction that $E=\bigcup_n F_n$ for some compact sets $F_n$. For each positive integer $n$, we will construct now a dimension function $h_n$ such that $\cH^{h_n}(F_n)=0$ and
\begin{equation}\label{eq:h_n_fast_decay}
    \lim_{\varepsilon\to0}\frac{\log h_n(\varepsilon)}{\log\varepsilon} = 0.
\end{equation}
Fix $n$ and denote $F=F_n$, where we dropped the dependence on $n$ to simplify the following notation. Since $F$ is a compact set with Hausdorff dimension zero, by definition, there is a finite open covering $\{U_{1,i}\}_i$ of $F$ such that $\sum_i |U_{1,i}|<1$. Let $r_1=\min_i|U_{1,i}|$ and $s_1=\max_i|U_{1,i}|$. Assuming $r_1,\dots,r_k,s_1,\dots,s_k$ are defined, by the same reason above, there is a finite open covering $\{U_{k+1,i}\}_i$ of $F$ such that 
\begin{equation}\label{eq:def_r_k}
    \sum_i |U_{k+1,i}|^{1/(k+1)}<r_k, \qquad\text{and}\qquad \max_i|U_{k+1,i}|<r_k^{1+1/k},
\end{equation}
so define $r_{k+1}=\min_i|U_{k+1,i}|$ and $s_{k+1}=\max_i|U_{k+1,i}|$. Finally, let $h_n$ be the right continuous non-decreasing function
\begin{equation*}
    h_n(\varepsilon)=\begin{cases}
        \varepsilon^{1/k}, & \text{if } r_k\leq\varepsilon < s_k; \\
        r_k^{1/k}, & \text{if } s_{k+1}\leq\varepsilon \leq r_k.
    \end{cases}
\end{equation*}

Note that $r_1\leq s_1<1$ since $\sum_i|U_{1,i}|<1$. Since \eqref{eq:def_r_k} implies $r_{k+1}\leq s_{k+1}<r_{k}^{1+1/k}$ for all $k\geq 1$, it follows by induction that $r_{k}\leq s_k<1$ for all $k\geq 1$ and also $r_{k+1}\leq r_1^{(1+1/1)\dotsb(1+1/k)}=r_1^{k+1}$. Since $r_1<1$, this shows that $r_k,s_k\to0$ as $k\to\infty$, so indeed $h_n$ is a non-decreasing function such that $\lim_{\varepsilon\to0}h_n(\varepsilon)=0$. The fact that $\cH^{h_n}(F)=0$ is a consequence of the construction: given any $\delta>0$, we let $k$ be sufficiently large so that $r_k^{1+1/k}<\delta$, so using \eqref{eq:def_r_k} with the covering $\{U_{k+1,i}\}_i$ of $F$ gives that
\begin{equation*}
    \sum_i h(\diam U_{k+1,i}) = \sum_i  |U_{k+1,i}|^{1/(k+1)}<r_k.
\end{equation*}
Since $r_k\to0$ as $k\to\infty$, we conclude that $\cH^{h_n}(F)=0$. Now to prove \eqref{eq:h_n_fast_decay}, just observe that for $s_{k+1}\leq\varepsilon<s_k$ one has that $0<\log h_n(\varepsilon)/\log\varepsilon\leq 1/k$, so \eqref{eq:h_n_fast_decay} follows from the fact that $s_k\to 0$ as $k\to\infty$.

Now that we have constructed dimension functions $h_n$ that satisfy \eqref{eq:h_n_fast_decay} and  $\cH^{h_n}(F_n)=0$, we will build another dimension function $h$ such that for each $n$, we have $h(\varepsilon)\leq h_n(\varepsilon)$ for all $\varepsilon>0$ sufficiently small but 
\begin{equation*}
    \lim_{\varepsilon\to0}\frac{\log h(\varepsilon)}{\log\varepsilon}=0.
\end{equation*}
In particular, it follows that $\cH^{h}(F_n)\leq\cH^{h_n}(F_n)=0$ for all $n$, but this is a contradiction to the hypothesis \eqref{eq:decay_ratio_general_topology} since $\cH^{h}(E)\leq\sum_n\cH^h(F_n)=0$.

To construct such a function, given any positive integer $n$, let $0<\delta_n<\delta_{n-1}$ be such that
\begin{equation*}
    \sup_{0<\varepsilon\leq\delta_n}\left\{\frac{\log h_1(\varepsilon)}{\log\varepsilon},\dots,\frac{\log h_n(\varepsilon)}{\log\varepsilon}\right\}<\frac{1}{n}.
\end{equation*}
Finally, define the right continuous dimension function $h(\varepsilon)=\min\{h_1(\varepsilon),\dots,h_n(\varepsilon)\}$ for $\delta_{n+1}\leq\varepsilon<\delta_n$. It is easy to check that $h$ satisfies the claimed properties, so the proof that $E$ is not $F_\sigma$ is done.

\end{proof}

\begin{lemma}\label{lem:topology}
The set $k^{-1}(3)$ is not $G_\delta$ and not $F_\sigma$, but it is $F_{\sigma\delta}$.
\end{lemma}

\begin{proof}
In general, for $t\in [0,\infty]$ such that $k^{-1}(t)$ is non-empty, the set $k^{-1}(t)$ is not a $G_\delta$ set, i.e., is not the intersection of countably many open sets. The reason why $k^{-1}(t)$ is not $G_\delta$ is because for any $s\in\R$, we have that
\begin{equation*}
    k^{-1}[s,\infty]=\bigcap_{n\geq 1}\bigcup_{q\geq 2}\bigcup_{p\in\Z}\left(\frac{p}{q}-\frac{1}{(s-1/n)q^2},\frac{p}{q}+\frac{1}{(s-1/n)q^2}\right)
\end{equation*}
is a dense $G_\delta$ set. If $k^{-1}(t)$ is $G_\delta$, then $k^{-1}(t)\cap k^{-1}[t+1,\infty]=\varnothing$ would be a dense $G_\delta$ set, being the intersection of two $G_\delta$ dense sets, a contradiction. 

Now we prove that the set $k^{-1}(3)$ is an $F_{\sigma\delta}$ set: it is the countable intersection of $F_{\sigma}$ sets where each $F_{\sigma}$ is the countable union of closed sets with empty interior. Indeed, since $k^{-1}[3+1/m,\infty]$ is a $G_\delta$ set, we have that $k^{-1}(-\infty,3+1/m)=\bigcup_{n\geq 1}X_{m,n}$ where each $X_{m,n}$ is a closed set with empty interior (since $k^{-1}(-\infty,t)$ contains no rational number for any $t\in\R$). On the other hand, the set $k^{-1}(-\infty,3)$ is actually countable, because $\{t<3:k^{-1}(t)\neq\varnothing\}=\{\ell_1=\sqrt{5}<\ell_2=2\sqrt{2}<\dots\}$ is countable and each $k^{-1}(\ell_n)$ is also countable because of Markov's theorem \cite[Theorem 16]{Bombieri}. In particular, since the complement of a point is the countable union of closed sets, the complement of the countable set $k^{-1}(-\infty,3)$ is an $F_{\sigma\delta}$, that is, it can be written as $k^{-1}[3,\infty]=\R\setminus k^{-1}(-\infty,3)=\bigcap_{m\geq 1}\bigcup_{n\geq 1} Y_{m,n}$ for some closed sets $Y_{m,n}$. Finally, we have that
\begin{equation*}
    k^{-1}(3)=k^{-1}(-\infty,3]\cap k^{-1}[3,\infty)=\bigcap_{m\geq 1}\bigcup_{n_1,n_2\geq 1}(X_{m,n_1}\cap Y_{m,n_2})
\end{equation*}
where each $X_{m,n_1}\cap Y_{m,n_2}$ is a closed set with empty interior.

Finally, notice that if $k^{-1}(3)$ is $F_\sigma$ then $k^{-1}(3)\cap[0,1]$ would be a bounded $F_\sigma$ set, but by \Cref{lem:infinite_case} we have $\cH^{h}(k^{-1}(3)\cap[0,1])=\infty$ (since $k^{-1}(3)$ is invariant by integer translations), which is a contradiction to \Cref{lem:general_topology}.

\end{proof}

The previous lemma contrasts with the situation for Liouville numbers, which form a dense $G_\delta$ set. Another major difference is that the set of Liouville numbers is invariant under translation by rational numbers. It is unsurprising that $k^{-1}(3)$ is not invariant by such translations, since adding a rational number to a continued fraction can drastically alter its coefficients. Mendès France \cite{MendesImprovement} proved, refining a result by Davenport \cite{Davenport}, that for any irrational $x$ and prime $p$, there is $0\leq a<p-1$ such that the continued fraction $x+\frac{a}{p}=[0;b_1,b_2,\dots]$ satisfies $b_n\geq p-1$ for a positive proportion of $n$, and consequently $k\left(x+\frac{a}{p}\right)\geq p-1$. As consequence, $k^{-1}(3)$ (and more generally $k^{-1}(t)\neq\varnothing$) can not be invariant by translations of the form $x\mapsto x+\frac{a}{b}$ with $\frac{a}{b}\in\Q\setminus\Z$, $\gcd(a,b)=1$, if $b$ is divisible by some prime $p\geq 5$.

To illustrate the previous phenomena with a concrete example, let $e_m = 6 m$ and consider the continued fraction
\begin{equation*}
    x=[0;1^{e_1},2,2,1^{e_2},2,2,1^{e_3},2,2,\dots].
\end{equation*}
By \Cref{lem:attainable_criteria} we have that $x\in k^{-1}(3)$. If we denote by $p_n/q_n$ the continued fraction of $x$, then it is an elementary exercise to show that for all $m\geq 1$ 
\begin{equation*}
    q_{e_1+\dots+e_{2m+1}+4m+1} \equiv 2 \pmod{4} \qquad\text{and}\qquad p_{e_1+\dots+e_{2m+1}+4m+1} \equiv 1 \pmod{4}.
\end{equation*}
In particular, for $n=e_1+\dots+e_{2m+1}+4m+1$, letting 
\begin{equation*}
    y=x+1/2, \qquad\text{and}\qquad (p,q)=((p_n+q_n/2)/2,q_n/2),
\end{equation*}
we will have that the inequality
\begin{equation*}
    \left|y-\frac{p}{q}\right|=\left|\left(x+\frac{1}{2}\right)-\left(\frac{p_n}{q_n}+\frac{1}{2}\right)\right|=\left|x-\frac{p_n}{q_n}\right|<\frac{1}{q_n^2}=\frac{1}{4q^2},
\end{equation*}
has infinitely many solutions $p/q\in\Q$. By the definition of the function $k$, we have that $k(x+1/2)=k(y)\geq 4$ (in fact $k(y)\geq 12$).

\section{Dimensionality of $K_3$}

For the Cantor set
\begin{multline*}
    K_3=\{[0;a_1,\dots,a_n,\dots]\ \mid\ \text{there exists $(a_{-n})_{n\geq 0}\in(\N_{>0})^{\N}$ such that} \\ 
[a_k;a_{k+1},\dots,]+[0;a_{k-1},a_{k-2},\dots]\leq 3,\forall k\in\Z\},
\end{multline*}
we obtain a reasonably explicit description of its generalized Hausdorff measures. However, due to additional combinatorial complications, the result is not as precise as in the case of $k^{-1}(3)$.

Recall that for any fixed $n\geq 1$, the set $\{I(w):w=(a_1,\dots,a_{n})\in\Sigma(3,n)\}$ is a covering of $K_3$ by closed intervals.

\begin{lemma}\label{lem:easy_case_K_3}
Let $h$ be a dimension function. If
\begin{equation*}
    \limsup_{\varepsilon\to 0}\frac{|\log h(\varepsilon)|}{\log|\log\varepsilon|}>3
\end{equation*}
then $\cH^h(K_3)=0$. 
\end{lemma}

\begin{proof}
By hypothesis, there is a positive $d>0$ and a decreasing sequence $\varepsilon_k$ such that $\varepsilon_k\to0$ when $k\to\infty$ such that
\begin{equation}\label{eq:fast_decay_h2}
    h(\varepsilon_k)\leq \frac{1}{|\log\varepsilon_k|^{3+d}}, \quad \text{for all $k$}.
\end{equation}

For a given $k$, define $n_k=\lceil\log_2 1/\varepsilon_k\rceil$. Taking $w=(a_1,\dots,a_{n_k})\in\Sigma(3,n_k)$, since $\diam I(w)<2^{-n_k}$, we have that
\begin{align*}
    \sum_{w\in\Sigma(3,n_k)}h(\diam I(w))&\leq |\Sigma(3,n_k)|\cdot h(\varepsilon_k) \\
    &\leq \frac{9n_k^3}{|\log\varepsilon_k|^{3+d}}\leq\frac{9n_k^3}{(\log2/2)^{3+d}n_k^{3+d}}=\frac{9}{(\log2/2)^{3+d}n_k^d}.
\end{align*}
Since $n_k\to\infty$ as $k\to\infty$, this proves the lemma.
\end{proof}

We can actually get a more precise estimation of the gauge functions that give zero measure. For example we can prove that the function $h(\varepsilon)=|\log\varepsilon|^{-3}$ gives $\cH^h(K_3)=0$.

The strategy for the proof is as follows. As we saw in the introduction, for a Christoffel pair $(\alpha,\beta)$ that corresponds to an interval $I(\alpha)$ of the covering of $K_3$, we can associate a unique rational number whose continued fraction $[0;b_1,\dots,b_m]$ encodes the construction of this pair. Similar to the proof that irrationals $[0;a_1,a_2,\dots]$ whose coefficients tend to $\infty$ have Hausdorff dimension 1/2, whenever a very large digit appears in the continued fraction of $[0;b_1,\dots,b_m]$, say $b_{j+1}\geq M+1$ for some appropriate chosen $M$, it is cheaper to use the alphabet $(\tilde{\alpha},\tilde{\beta})$ associated to $[0;b_1,\dots,b_j,M]$ and cover the interval $I(\alpha)$ with $I(\tilde{\alpha})$ instead. 

\begin{lemma}\label{lem:upper_bound_K_3}
For the dimension function 
\begin{equation*}
    g(\varepsilon)=\frac{\log|\log\varepsilon|}{|\log\varepsilon|^3(\log\log|\log\varepsilon|)},
\end{equation*}
we have that $\cH^{g}(K_3)<\infty$.
\end{lemma}

\begin{proof}
We will consider initially the covering of $K_3$ given by the intervals $I(a_1,\dots,a_N)$ where $(a_1,\dots,a_N)\in\Sigma(3,N)$. Since some of these intervals are very small, we will show that they can be covered more efficiently by bigger intervals $I(a_1,\dots,a_\ell)$ associated with shorter prefixes.

Given any word $(a_1,\dots,a_N)\in\Sigma(3,N)$, we know by \Cref{lem:3w} that there is a word $w\in P$ with $N\leq |w|\leq 3N$ such that $(a_1,\dots,a_N)$ is a subword of $w$. In particular, if $w\neq b$, we can write $w=\alpha$ for some $(\alpha,\beta)\in\overline{P}$ with $|\alpha|>|\beta|$. Writing $w=w_1\dots w_{|w|}\in\{1,2\}^{|w|}$, since $\alpha\alpha$ is a subfactor of $\alpha^\infty\in\Sigma_3$, all cyclic permutations $w_i\dots w_{|w|}w_1\dots w_{i-1}$, $1\leq i\leq |w|$ appear in $\Sigma(3,|w|)$ and thus appear as prefixes of continued fractions in $K_3$. In particular we have that $a_1\dots a_N$ is a prefix of $w_i\dots w_{|w|}w_1\dots w_{i-1}$ for some $1\leq i\leq |w|$. We will extend each $(a_1,\dots,a_N)$ so that it coincides with that permutation $(a_1,\dots,a_{N^\prime})=(w_i,\dots,w_{|w|},w_1,\dots,w_{i-1})$ for some $N^\prime\leq 3N$. We will refine further this covering by considering all $\alpha\in P$ with $N\leq N^\prime=|\alpha| \leq 6N$. 

On the other hand, as explained in \Cref{subsec:Sigma_3}, there is a unique continued fraction $[0;b_1,\dots,b_m]$ that codifies the Christoffel pair $(\alpha,\beta)\in\overline{P}$ and satisfies $q_m=|\alpha|/2$. Notice that the denominator of such continued fraction is between $N/2$ and $3N$. Let $M$ be a large positive value that depends on $N$ that we will choose later. We have two different situations, either there is a coefficient $b_{j+1}>M$ with $j$ minimal or $b_j\leq M$ for all $1\leq j\leq m$.

Now we will consider the first situation. Let us count how many such continued fractions with denominator between $N/2$ and $3N$ and such that $M<b_{j+1} \leq 2M$ is the first coefficient larger than $M$ exist. Writing $[0;b_1,\dots,b_j,M]=\frac{p_r}{q_r}$, any such fraction is of the form $[0;b_1,\dots,b_j,M+s,\frac{q}{p}]$ where $0<s\leq M$ and $0\leq p/q\leq 1$. Writing $[0;b_1,\dots,b_j,M+s]=\frac{p_r+sp_{r-1}}{q_r+sq_{r-1}}$, we have that 
\begin{equation*}
    [0;b_1,\dots,b_j,M+s,q/p]=\frac{q(p_r+sp_{r-1})+pp_{r-1}}{q(q_r+sq_{r-1})+pq_{r-1}}.
\end{equation*}
Since $q(q_r+sq_{r-1})+pq_{r-1}\leq q(q_r+Mq_{r-1})+qq_{r-1}<3qq_r$, we see that any choice of $\frac{N}{2q_r}<q<\frac{N}{q_r}$, $0\leq p\leq q$ and  $0< s\leq M$, gives $q(q_r+sq_{r-1})+pq_{r-1}\in[N/2,3N]$. Similarly we have $\frac{N}{6q_r}<q<\frac{3N}{q_r}$. Therefore we have at least $(\frac{N}{2q_r})^2$ and at most $\left(\frac{3N}{q_r}\right)^2$ choices for $0\leq p\leq q$, and we have $M$ choices for $s$. In summary, we have between $\frac{N^2M}{4q_r^2}$ and $\frac{9N^2M}{q_r^2}$ such continued fractions $[0;b_1,\dots,b_m]=[0;b_1,\dots,b_j,M+s,\frac{q}{p}]$. Now, since permutations of $\alpha$ also give different possibilities for $(a_1,\dots,a_{N^\prime})$, we have between $\frac{N^3M}{4q_r^2}$ and $\frac{54N^3M}{q_r^2}$ intervals $I(a_1,\dots,a_{N^\prime})$ with $(a_1,\dots,a_{N^\prime})$ a cyclic permutation of $\alpha$ associated to a unique $[0;b_1,\dots,b_j,M+s,q/p]$. All of them have diameter at most $\left(\frac{3+\sqrt{5}}{2}\right)^{-6N+1}$. On the other hand we know that $|\Sigma(3,N^\prime)|\leq 9(N^\prime)^3\leq 1944N^3$, so if we denote
\begin{equation*}
    \mathcal{C}=\left\{(p_r,q_r):\substack{\displaystyle 1\leq p_r\leq q_r\leq 3N, \gcd(p_r,q_r)=1, \\ \displaystyle p_r/q_r=[0;b_1,\dots,b_j,M],1\leq b_i\leq M, 1\leq i\leq j}\right\},
\end{equation*}
we have summing over all such intervals
\begin{equation}\label{eq:36M}
    \sum_{(p_r,q_r)\in\mathcal{C}}\frac{1}{q_{r}^2} \leq \frac{|\Sigma(3,6N)|}{54N^3M}\leq\frac{36}{M}.
\end{equation}

Now we will change all these intervals by a more effective covering. Let us denote $(\tilde{\alpha},\tilde{\beta})\in\overline{P}$ the pair corresponding to $[0;b_1,\dots,b_j,M]=\frac{p_r}{q_r}$. In particular we have that either $\tilde{\alpha}=\hat{\alpha}\tilde{\beta}^M$ with $q_r=|\tilde{\alpha}|/2$ or $\tilde{\beta}=\tilde{\alpha}^M\hat{\beta}$ with $q_r=|\tilde{\beta}|/2$ for some $(\hat{\alpha},\hat{\beta})\in\overline{P}$. Let us denote $(\alpha^\prime,\beta^\prime)\in\{(\hat{\alpha}\tilde{\beta}^M,\tilde{\beta}),(\tilde{\alpha},\tilde{\alpha}^M\hat{\beta})\}$ accordingly. Notice that $\alpha$ can be written in the alphabet $(\alpha^\prime,\beta^\prime)$ and that moreover it begins with $\alpha^\prime$ and ends with $\beta^\prime$. First, suppose that $\min\{|\hat{\alpha}|,|\hat{\beta}|\}>2$. Since $(a_1,\dots,a_{N^\prime})$ is a permutation of $\alpha$, we have that there is a prefix of $(a_1,\dots,a_{N^\prime})$ that is subword of length $2q_r$ of $\tilde{\theta}\tilde{\theta}$ where $\tilde{\theta}\in\{\hat{\alpha}\tilde{\beta}^M,\tilde{\alpha}^M\hat{\beta}\}$. Hence, there is a prefix $(a_1,\dots,a_\ell)$ of $(a_1,\dots,a_{N^\prime})$ that is a permutation of $\tilde{\theta}$. Each of these intervals $I(a_1,\dots,a_\ell)$ has diameter at most $\left(\frac{3+\sqrt{5}}{2}\right)^{-2q_r+1}\leq\left(\frac{3+\sqrt{5}}{2}\right)^{-2M+1}$ by \eqref{eq:length_of_intervals_12} and we have precisely $2q_r$ permutations. Therefore their $g$--measure is at most
\begin{equation}\label{eq:bound_small_intervals}
    \sum_{(p_r,q_r)\in\mathcal{C}}2q_rg\left(\left(\frac{3+\sqrt{5}}{2}\right)^{-2q_r+1}\right)\leq C^\prime\sum_{(p_r,q_r)\in\mathcal{F}}\frac{1}{q_r^2}\frac{\log q_r}{\log\log q_r},
\end{equation}
for some constant $C^\prime>0$. Now we will choose $M$ to be equal to $M=C(\log N)/(\log\log N)$ for some absolute constant $C>0$ we will choose later. In particular, using \eqref{eq:36M} and \eqref{eq:bound_small_intervals}
\begin{equation*}
    \sum_{(p_r,q_r)\in\mathcal{C}}2q_rg\left(\left(\frac{3+\sqrt{5}}{2}\right)^{-2 q_r+1}\right)\leq\frac{C^\prime}{C}\sum_{(p_r,q_r)\in\mathcal{F}}\frac{M}{q_r^2}=\frac{36C^\prime}{C}<\infty.
\end{equation*}

If $\min\{|\hat{\alpha}|,|\hat{\beta}|\}=2$, then $j+1=1$ and $(\hat{\alpha},\hat{\beta})=(a,b)$ and $(a_1,\dots,a_{N^\prime})$ has no prefix that is permutation of $\tilde{\theta}\tilde{\theta}$, it is because either $2b^{M+1}$, $b^{M+1}$ or $a^{M+1}$ is a prefix of $(a_1,\dots,a_{N^\prime})$. Since the diameters of the intervals $I(2b^{M+1})$, $I(b^{M+1})$, $I(a^{M+1})$ go to zero provided $M$ goes to infinity as $N$ goes to infinity, the contribution of these intervals can be ignored. In conclusion the contribution of all the intervals in the first situation is bounded by some absolute constant.

On the other hand, in the second situation, if all coefficients are $b_j\leq M$ for $1\leq j\leq m$, then we will not change $I(a_1,\dots,a_{N^\prime})$. By \eqref{eq:Q_n(x)} we have at most $2(3N)^{2\dim_H(C(M))}$ continued fractions $[0;b_1,\dots,b_m]$ with denominator at most $3N$ and $1\leq b_i\leq M$ for $1\leq i\leq m$. For each such continued fraction, it corresponds a unique $\alpha\in P$ and at most $|\alpha|\leq 6N$ permutations. Using Jarnik's estimate there is an absolute constant $C_1>0$ such that $\dim_H(C(M))<1-C_1/M$, therefore we have at most $4(3N)^{2\dim_H(C(M))+1}\leq 108 N^3N^{-2C_1/M}$ intervals. Now we will fix $C=2C_1$, so we have at most $108N^3/\log N$ such intervals with diameter at most $\left(\frac{3+\sqrt{5}}{2}\right)^{-6N+1}$. Therefore, they have $g$--measure at most
\begin{equation*}
    \frac{108N^3}{\log N}g\left(\left(\frac{3+\sqrt{5}}{2}\right)^{-6N+1}\right)\leq\hat{C}\frac{1}{\log\log N},
\end{equation*}
for some absolute constant $\hat{C}>0$. This shows that the contribution of these intervals goes to zero.

\end{proof}

\begin{lemma}\label{lem:lower_bound_K_3}
For a large enough constant $\hat{C}>0$, the dimension function 
\begin{equation*}
    g(\varepsilon)=\frac{\exp(\hat{C}\sqrt{\log|\log\varepsilon|})}{|\log\varepsilon|^3}
\end{equation*}
satisfies $\cH^g(K_3)=\infty$. In particular for any dimension function $h$ such that
\begin{equation*}
    \limsup_{\varepsilon\to 0}\frac{|\log h(\varepsilon)|}{\log|\log\varepsilon|}<3,
\end{equation*}
we have $\cH^h(K_3)=\infty$. 
\end{lemma}

\begin{proof}
We will use \Cref{cantor}. We choose $C_n=\bigcup\{I(w):w=(a_1,\dots,a_{n})\in\Sigma(3,2^{n^2})\}$ to be our coverings so in particular $K_3=\bigcap_{n\geq 1}C_n$. We claim that there is a constant $C>0$ such that any interval $I(w)$ with $w\in\Sigma(3,2^{n^2})$ is subdivided in at least $C\cdot 2^{3\cdot(2n+1)}$ intervals $I(\tilde{w})$ with $\tilde{w}\in\Sigma(3,2^{(n+1)^2})$, or equivalently, each $w\in\Sigma(3,2^{n^2})$ has at least $C\cdot 2^{3(2n+1)}$ different such continuations $\tilde{w}\in\Sigma(3,2^{(n+1)^2})$.  

Now we will prove the claim. Given any word $w\in\Sigma(3,n^\prime)$, by \Cref{lem:3w} there is a Christoffel word $\alpha\in P\cup\{a,b\}$ with $2m=|\alpha|<3n^\prime$ such that $w$ is subfactor of $\alpha$. By \Cref{lem:farey_neighbors}, there is another Christoffel word $\beta\in P\cup\{a,b\}$ such that $(\alpha,\beta)\in\overline{P}$ is an alphabet and such that $\Theta(\alpha)<\Theta(\beta)$ are Farey neighbors in $\mathcal{F}_m$. In particular, $\max\{|\alpha|_{a,b},|\beta|_{a,b}\}\leq m<|\alpha|_{a,b}+|\beta|_{a,b}$. 

Denote by $L_k$ the set of all words $l_1\dots l_k\in\{a,b\}^k$ that appear as subfactors of some word in $\Sigma_3$. It is known \cite{Reutenauer2006} that $L_k$ is the same as the Sturmian factors of length $k$ and it is known \cite{countingSturmianwords} that we have $|L_k|\sim k^3/\pi^2$. For our purposes, it suffices to prove that there is a constant $\tilde{C}>0$ such that $|L_k|\geq 2\tilde{C}k^3$ for all $k$. Let us give an alternative proof of this fact that it is useful for the next section. Recall that given a Christoffel word $w\in P\cup\{a,b\}$ we can associate a unique Farey fraction $\Theta(w)=\frac{|w|_b}{|w|_{a,b}}$ where $|w|_b$ denotes the number of letters $b$ in $w$ and $|w|_{a,b}$ is the length of $w$ over the alphabet $\{a,b\}$. In particular if $\mathcal{F}_k$ denotes the $k$--term of the Farey sequence, we can find at least $|\mathcal{F}_k|\sim \frac{3k^2}{\pi}$ elements in $L_k$. Since cyclic permutations also appear as subfactors of longer Christoffel words and since $|\mathcal{F}_k|=|\mathcal{F}_{k-1}|+\varphi(k)$ where $\varphi$ denotes the Euler function, we have the known formula $|L_k|=1+\sum_{i=1}^k(k-i+1)\varphi(i)\sim\frac{k^3}{\pi^2}$.

On the other hand, since the operation of exchanging $a$ with $b$ preserves factors in $L_k$, precisely half of the words in $L_k$ begin with $a$, say $\{\eta_1,\dots,\eta_{|L_k|/2}\}$.

Since $(\alpha,\beta)\in\overline{P}$, there are some inner renormalization operators $R_1,\dots,R_\ell\in\{U,V\}$ such that $\alpha=R(a)$ and $\beta=R(b)$ where $R=R_1\dotsb R_\ell$. In particular, each $R(\eta_j)$, $1\leq j\leq |L_k|/2$ is also a Sturmian factor of length over the alphabet $\{1,2\}$ of at most $mk<\frac{3}{2}n^\prime k$. Each such factor $R(\eta_j)$ begins with $\alpha$ and since $w$ is a subfactor of $\alpha$, each $R(\eta_j)$ provides a continuation of $w$ of length over the alphabet $\{1,2\}$ at most $\frac{3}{2}n^\prime k$. Since any Sturmian factor can be extended arbitrarily, we conclude that any $w\in\Sigma(3,n^\prime)$ has at least $|L_k|/2\geq\tilde{C
}k^3$ continuations $\tilde{w}\in\Sigma(3,3n^\prime k)$. Consequently, to finish the proof of the claim, we choose $k=\lfloor 2^{2n+2} /3\rfloor$ and $n^\prime=2^{n^2}$.

Notice that by \eqref{eq:length_of_intervals_12} the size of each interval of $C_n$ is at least $\varepsilon_n = (3+2\sqrt{2})^{-2^{n^2}-1}$. According to the notation of \Cref{cantor}, we will have $N_n=\prod_{j<n}(C\cdot 2^{3(2j+1)})=C^{n-1}2^{3n^2-3}$. In particular, for a large enough constant $\hat{C}>0$, the dimension function 
\begin{equation*}
    g(\varepsilon)=\frac{\exp(\hat{C}\sqrt{\log|\log\varepsilon|})}{|\log\varepsilon|^3}
\end{equation*}
satisfies 
\begin{equation*}
    \lim_{n\to\infty}N_n\cdot g(\varepsilon_{n+1})\geq\lim_{n\to\infty}\frac{\exp(n\log(C/64)+\hat{C}(n+1)\sqrt{\log2})}{64C\log(3+2\sqrt{2})^3}=\infty,
\end{equation*}
so by \Cref{cantor} we obtain $\mathcal{H}^g(K_3)=\infty$. Moreover, the function $g$ clearly satisfies
\begin{equation*}
    \lim_{\varepsilon\to 0}\frac{|\log g(\varepsilon)|}{\log|\log\varepsilon|}=\lim_{\varepsilon\to 0}\left\vert\frac{\hat{C}}{\sqrt{\log|\log\varepsilon|}}-3\right\vert=3.
\end{equation*}
Finally, by the hypothesis on the dimension function $h$, for all sufficiently small $\varepsilon>0$ we have $h(\varepsilon)>g(\varepsilon)$. Consequently $\mathcal{H}^h(K_3)\geq\mathcal{H}^g(K_3)=\infty$.
\end{proof}

\section{Dimensionality of $B_3$}

The following is a consequence of the main theorem of \cite{cao2026badlyapproximablenumbers}. Recall that given a word $w=w_1\dots w_n\in\{a,b\}^n$, we denote $w^{+}=w_2\dots w_n$ and $w^T=w_n\dots w_1$. 

\begin{theorem}\label{thm:full_characterization}
Let $x=[x_0;x_1,x_2,\dots]$ be such that $\left| x-\frac{p}{q} \right| < \frac{1}{3q^2}$ has only finitely many solutions and let $N\in\N=\{0,1,\dots\}$ minimal such that $\lambda_n(x)\leq3$ for all $n\geq N+1$ and $k(x)=3$. Then there is a sequence of alphabets $(\alpha_{n+1},\beta_{n+1})\in \{\overline{U}(\alpha_{n},\beta_{n}),\overline{V}(\alpha_{n},\beta_{n})\}$ with both renormalization operators $\overline{U}$ and $\overline{V}$ appearing infinitely many times and such that:
\begin{itemize}
    \item if $x_{N+1}=1$ then $x_{N+1}x_{N+2}x_{N+3}\ldots=\lim_{n\to\infty}\beta_n^T$ and if $N\geq 1$ then $[0;x_N,\dots,x_1]<[0;\lim_{n\to\infty}\alpha_n]$; or
    \item if $x_{N+1}=2$ then $x_{N+1}x_{N+2}x_{N+3}\ldots=\lim_{n\to\infty}2\alpha_n^{+}$ and if $N\geq1$ then $[0;x_N,\dots,x_1]<[0;2,\lim_{n\to\infty}\beta_n^T]$.
\end{itemize}
Conversely, for any such continued fraction $x=[x_0;x_1,x_2,\dots]$ the inequality $\left| x-\frac{p}{q} \right| < \frac{1}{3q^2}$ has at most $N$ solutions.
\end{theorem}

Let us denote by $B_3(N)$ the set of $x\in B_3\cap(0,1)$ such that $N\in\N$ is minimal with $\lambda_n(x)\leq3$ for all $n\geq N+1$. Since 
\begin{equation*}
    B_3(N)\subset\bigcup_{(x_1,\dots,x_N)\in(\N_{>0})^N}\{[0;x_1,\dots,x_N,1/\gamma]:\gamma\in B_3(0)\cap(0,1)\},
\end{equation*}
we have that $\cH^h(B_3(N))\leq\cH^h(B_3(0))$. On the other hand, we clearly have $\cH^h(B_3)\geq\cH^h(B_3(0))$. Therefore, to estimate the generalized Hausdorff measure of $B_3$, we will consider instead $B_3(0)$, that is, the set
\begin{equation*}
    B_3(0) = \left\{x\in\R\setminus\Q:\left|x-\frac{p}{q}\right|\geq\frac{1}{3q^2} \text{  for all }\frac{p}{q}\in\Q\right\}.
\end{equation*}

Now another reduction can be made to the study of the dimensionality of the set $B_3(0)$. Indeed, given a sequence $d_1,d_2,\dots$ of positive integers, denote $i_n=d_n$ and $\overline{R}_n=\overline{U}$ if $n$ is even and $i_n=d_n-1$ and $\overline{R}_n=\overline{V}$ is $n$ is odd, then we have that $(x_0,x_1,x_2,\dots)=\lim_{n\to\infty}\alpha_n$   where $(\alpha_n,\beta_n)=\left(\overline{R_n}^{i_n}\dotsb\overline{R_1}^{i_1}\right)(a,b)$ defines uniquely the number $[0;x_1,x_2,\dots]$ of $B_3(0)$. It turns out that if we change $(x_0,x_1)=(2,2)$ by $(1,1)$, then $(1,1,x_2,x_3,\dots)=\lim_{n\to\infty}\beta_n^T$, so $[0;1,1,x_2,x_3,\dots]$ also defines the corresponding number of $B_3(0)$. To see this, denote $\alpha^b=b\alpha^{+}$, so we have that $(1,1,x_2,\dots)$ begins with $\alpha_{2n}^b=(\alpha_{2n-1}^{d_{2n-1}}\beta_{2n-1})^b(\alpha_{2n-1}^{d_{2n-1}-1}\beta_{2n-1})^{d_{2n}-1}=(\beta_{2n-1}^T(\alpha_{2n-1}^T)^{d_{2n-1}-1}\alpha_{2n-1}^b)(\alpha_{2n-1}^{d_{2n-1}-1}\beta_{2n-1})^{d_{2n}-1}$, begins with $\beta_{2n-1}^T$. Since this holds for all $n$, the sequence $(1,1,x_2,\dots)$ coincides with the limit $\lim_{n\to\infty}\beta_n^T$. 

Therefore, it suffices to consider numbers of the form $[0;2,\lim_{n\to\infty}\alpha_n^{+}]$ in $B_3(0)$. Moreover, since $\frac{1}{2+B_3(0)}$ is a Bi-Lipschitz copy of $B_3(0)$, it suffices to consider the numbers of the form $[0;\lim_{n\to\infty}\alpha_n]$ where $\alpha_n$ is obtained through a sequence of renormalization operators where both $\overline{U}$ and $\overline{V}$ appear infinitely many times. 

The same proof of \Cref{lem:upper_bound_K_3} but ignoring the permutation of factors give the following. 
\begin{lemma}
For the dimension function 
\begin{equation*}
    g(\varepsilon)=\frac{\log|\log\varepsilon|}{|\log\varepsilon|^2(\log\log|\log\varepsilon|)},
\end{equation*}
we have that $\cH^{g}(B_3(0))<\infty$. In particular for any dimension function $h$ such that
\begin{equation*}
    \liminf_{\varepsilon\to 0}\frac{|\log h(\varepsilon)|}{\log|\log\varepsilon|}>2,
\end{equation*}
we have $\cH^h(B_3)=0$. 
\end{lemma}

It is clear that the statement of the previous theorem does not imply the second part of \Cref{thm:very_badly_approximable_case}, so one instead argues as in \Cref{lem:easy_case_K_3}. 

\begin{lemma}
    Let $h$ be a dimension function. If
\begin{equation*}
    \limsup_{\varepsilon\to 0}\frac{|\log h(\varepsilon)|}{\log|\log\varepsilon|}>2
\end{equation*}
then $\cH^h(B_3)=0$. 
\end{lemma}

\begin{proof}
    Recall that we already reduced the problem to consider numbers of the form $[0;\lim_{n\to\infty}\alpha_n]$ where $\alpha_n$ is obtained through a sequence of renormalization operators where both $\overline{U}$ and $\overline{V}$ appear infinitely many times. We can cover these numbers with the intervals $\{I(\alpha):m\leq |\alpha|< 2m,\alpha\in P\}$ where $P$ is the set of lower Christoffel words. Except for the constant words $a\dots a$ and $b\dots b$, by \Cref{lem:farey_neighbors} we can write $\alpha=uv$ where $(u,v)\in\overline{P}$ where $2|uv|\in\mathcal{F}_{2|uv|}$, the $2|uv|$--term of the Farey sequence. Since $|\mathcal{F}_{2m}|\sim \frac{12m^2}{\pi}$, we can cover our target set by roughly $\frac{12m^2}{\pi}$ intervals $I(\alpha)$ all of them with sizes at most $\left(\frac{3+\sqrt{5}}{2}\right)^{-2m+1}$ by \eqref{eq:length_of_intervals_12}. The rest of the proof proceeds as in \Cref{lem:easy_case_K_3}. 
\end{proof}

Finally, the same proof of \Cref{lem:lower_bound_K_3} but ignoring the permutation of factors give the following. 

\begin{lemma}
For a large enough constant $\tilde{C}>0$, the dimension function 
\begin{equation*}
    g(\varepsilon)=\frac{\exp(\tilde{C}\sqrt{\log|\log\varepsilon|})}{|\log\varepsilon|^2}
\end{equation*}
satisfies $\cH^g(B_3(0))=\infty$. In particular for any dimension function $h$ such that
\begin{equation*}
    \limsup_{\varepsilon\to 0}\frac{|\log h(\varepsilon)|}{\log|\log\varepsilon|}<2,
\end{equation*}
we have $\cH^h(B_3)=\infty$. 
\end{lemma}

\section*{Acknowledgements}

Nicolas Angelini gratefully acknowledges the Instituto de Matemática Pura e Aplicada (IMPA) for hosting his research visit, and the  Academia de Ciencias de América Latina (ACAL) for financial support during this stay. We also thank the referee for their thorough and careful reading of the paper.

\section*{Funding}
The research of the first author was partially supported by CNPq and FAPERJ. The resarch of the second author was partially supported by CAPES and FAPERJ.
The research of the third author was supported by \textit{Programa de Subvenciones para Pasantías Cortas de Investigadores Jóvenes} of the Academia de Ciencias de América Latina (ACAL), and was also partially supported by grants PICT 2022-4875 (ANPCyT), PIP 202287/22 (CONICET), and PROICO 3-0720 ``An\'alisis Real y Funcional. Ec. Diferenciales''.

\bibliographystyle{plain}
\bibliography{bibliography}

\begin{thebibliography}{10}

\bibitem{MassTransferencePrinciple}
Victor Beresnevich and Sanju Velani.
\newblock A mass transference principle and the {D}uffin-{S}chaeffer conjecture
  for {H}ausdorff measures.
\newblock {\em Ann. of Math. (2)}, 164(3):971--992, 2006.

\bibitem{Bombieri}
Enrico Bombieri.
\newblock Continued fractions and the {M}arkoff tree.
\newblock {\em Expo. Math.}, 25(3):187--213, 2007.

\bibitem{cao2026badlyapproximablenumbers}
Zhe Cao, Harold Erazo, and Carlos~Gustavo Moreira.
\newblock On very badly approximable numbers, 2026.
\newblock ArXiv:2602.09700.

\bibitem{CF}
Thomas~W. Cusick and Mary~E. Flahive.
\newblock {\em The {M}arkoff and {L}agrange spectra}, volume~30 of {\em
  Mathematical Surveys and Monographs}.
\newblock American Mathematical Society, Providence, RI, 1989.

\bibitem{DFSU}
Tushar Das, Lior Fishman, David Simmons, and Mariusz Urba{\'n}ski.
\newblock Hausdorff dimensions of perturbations of a conformal iterated
  function system via thermodynamic formalism: T. das et al.
\newblock {\em Selecta Mathematica}, 29(2):19, 2023.

\bibitem{Davenport}
H.~Davenport.
\newblock A remark on continued fractions.
\newblock {\em Michigan Math. J.}, 11:343--344, 1964.

\bibitem{NoBorelMeasure}
M\'arton Elekes and Tam\'as Keleti.
\newblock Borel sets which are null or non-{$\sigma$}-finite for every
  translation invariant measure.
\newblock {\em Adv. Math.}, 201(1):102--115, 2006.

\bibitem{EGRS2024}
Harold Erazo, Rodolfo Guti\'errez-Romo, Carlos~Gustavo Moreira, and Sergio
  Roma\~na.
\newblock Fractal dimensions of the {M}arkov and {L}agrange spectra near 3.
\newblock {\em J. Eur. Math. Soc. (JEMS)}, 28(9):3983--4046, 2026.

\bibitem{Fogg}
N.~Pytheas Fogg.
\newblock {\em Substitutions in dynamics, arithmetics and combinatorics},
  volume 1794 of {\em Lecture Notes in Mathematics}.
\newblock Springer, Berlin, 2002.

\bibitem{Hensley2}
Doug Hensley.
\newblock Continued fraction {C}antor sets, {H}ausdorff dimension, and
  functional analysis.
\newblock {\em J. Number Theory}, 40(3):336--358, 1992.

\bibitem{Jarnik1928}
Vojt{\v{e}}ch Jarn{\'i}k.
\newblock Zur metrischen theorie der diophantischen approximationen.
\newblock {\em Prace Matematyczno-Fizyczne}, 36(1):91--106, 1928/1929.

\bibitem{LMMSbook}
Davi Lima, Carlos Matheus, Carlos~G. Moreira, and Sergio Roma\~{n}a.
\newblock {\em Classical and dynamical {M}arkov and {L}agrange
  spectra---dynamical, fractal and arithmetic aspects}.
\newblock World Scientific Publishing Co. Pte. Ltd., Hackensack, NJ, 2021.

\bibitem{Mauduit}
Christian Mauduit and Carlos~Gustavo Moreira.
\newblock Generalized {H}ausdorff dimensions of sets of real numbers with zero
  entropy expansion.
\newblock {\em Ergodic Theory Dynam. Systems}, 32(3):1073--1089, 2012.

\bibitem{MendesImprovement}
Michel Mend\`es~France.
\newblock On a theorem of {D}avenport concerning continued fractions.
\newblock {\em Mathematika}, 23(2):136--141, 1976.

\bibitem{countingSturmianwords}
Filippo Mignosi.
\newblock On the number of factors of {S}turmian words.
\newblock {\em Theoret. Comput. Sci.}, 82(1):71--84, 1991.

\bibitem{M:geometric_properties_Markov_Lagrange}
Carlos~Gustavo Moreira.
\newblock Geometric properties of the {M}arkov and {L}agrange spectra.
\newblock {\em Ann. of Math. (2)}, 188(1):145--170, 2018.

\bibitem{Ol1}
L.~Olsen.
\newblock On the exact {H}ausdorff dimension of the set of {L}iouville numbers.
\newblock {\em Manuscripta Math.}, 116(2):157--172, 2005.

\bibitem{Ol2}
L.~Olsen and Dave~L. Renfro.
\newblock On the exact {H}ausdorff dimension of the set of {L}iouville numbers.
  {II}.
\newblock {\em Manuscripta Math.}, 119(2):217--224, 2006.

\bibitem{PalisTakens}
Jacob Palis and Floris Takens.
\newblock {\em Hyperbolicity and sensitive chaotic dynamics at homoclinic
  bifurcations}, volume~35 of {\em Cambridge Studies in Advanced Mathematics}.
\newblock Cambridge University Press, Cambridge, 1993.
\newblock Fractal dimensions and infinitely many attractors.

\bibitem{Reutenauer2006}
Christophe Reutenauer.
\newblock On {M}arkoff's property and {S}turmian words.
\newblock {\em Math. Ann.}, 336(1):1--12, 2006.

\bibitem{Reutenauerbook}
Christophe Reutenauer.
\newblock {\em From {C}hristoffel words to {M}arkoff numbers}.
\newblock Oxford University Press, Oxford, 2019.

\bibitem{Rogers}
C.~A. Rogers.
\newblock {\em Hausdorff measures}.
\newblock Cambridge Mathematical Library. Cambridge University Press,
  Cambridge, 1998.
\newblock Reprint of the 1970 original, With a foreword by K. J. Falconer.

\bibitem{GuguVillamil}
Christian~C. Silva~V. and Carlos~G. Moreira.
\newblock Hausdorff dimension of some subsets of the lagrange and markov
  spectra near 3, 2025.
\newblock ArXiv:2504.20300.

\end{thebibliography}

\end{document}